\documentclass{amsart}[11pt]
\usepackage{myAMSpackages}
\usepackage{Universal}
\usepackage{MyAMSEnvironments}
\usepackage{citesort}

\newcommand{\twoh}{\ensuremath{I\negmedspace I_0}\xspace}

\newcommand{\hess}[2][{}]{\ensuremath{{\mathcal{X}_{#1}^{#2}}}\xspace}
\newcommand{\wh}[1]{\ensuremath{{\widehat{#1}}}\xspace}
\newcommand{\Lie}[1]{\ensuremath{\mathcal{L}_{#1}}\xspace}

\newcommand{\divE}{\ensuremath{\text{div}_\Sigma}\xspace}
\newcommand{\omj}[1]{\omega^{(#1)}}
\newcommand{\etj}[1]{\eta^{(#1)}}
\newcommand{\mb}[1]{{\bar{#1}}}
\newcommand{\Cb}[2][{\ub}]{C^{#1}_{#2}}
\newcommand{\Gm}{\Gamma}
\newcommand{\om}{{\omega}}
\newcommand{\junk}{\cdot \text{junk}}
\newcommand{\F}{{\mathcal{F}}}
\newcommand{\GF}[1][{}]{{\mathcal{GF}_{#1}}}
\newcommand{\G}{{\mathcal{G}}}
\newcommand{\HF}[1][{}]{{\mathcal{GF}_{V#1}}}
\newcommand{\sM}{\s(M)}
\newcommand{\sMo}{\s_0(M)}
\newcommand{\s}{{\mathcal{S}}}
\newcommand{\vf}{\rho}
\newcommand{\vfo}{{\rho_0}}
\newcommand{\vfc}{\wh{\rho}}
\newcommand{\vfoc}{\wh{\rho}_0}

\newcommand{\pr}[1]{\ensuremath{#1^{\prime}}\xspace}
\newtheorem{MThm}{Theorem}

\newtheorem{Thm}{Theorem}

\begin{document}

\title{Variation of Perimeter Measure in sub-Riemannian geometry.}
\author{Robert K. Hladky}
\address{University of Rochester, Rochester, NY 14627}
\email{hladky@math.rochester.edu}
\author{Scott D. Pauls}
\address{Dartmouth College, Hanover, NH 03755}
\email{scott.pauls@dartmouth.edu}
\thanks{Both authors are partially supported by NSF grant DMS-0306752}

\begin{abstract}
We derive a formula for the first variation of horizontal perimeter measure for $C^2$ hypersurfaces of completely general sub-Riemannian manifolds, allowing for the existence of characteristic points. For $C^2$ hypersurfaces in vertically rigid sub-Riemannian manifolds we also produce a second variation formula for variations supported away from the characteristic locus. This variation formula is used to show the bubble sets in \hn{2} are stable under volume preserving variations. 
\end{abstract}

\maketitle

\section{Introduction}
Optimization problems lie at the heart of many pure and applied
problems, two of which, the minimal and isoperimetric surface
problems, have played a central role in mathematical development over
the last century.  In the last decade, there has been increasing
interest in these problems in the setting of sub-Riemannian spaces.
This interest is driven in part by applications of
sub-Riemannian geometry to optimal control as well as to
more novel applications such as the recent sub-Riemannian model of the
primary visual cortex \cite{V:Hoffman,V:PT,V:P,CS}.  We are also motivated by a deep
conjecture of Pansu \cite{Pansu} concerning the isoperimetric profile
in the sub-Riemannian Heisenberg group, which has seen a great deal of
recent activity with many partial results \cite{DGN,LR,RR,MR} (see
also \cite{CDPT} for an overview of this problem).  

One of the basic approaches to such optimization problems uses the
tools of the calculus of variations to determine the geometric and
analytical properties of their solutions.  Recent investigations of minimal and isoperimetric surfaces have focused on
this approach with many authors deriving first and second variation
formulae.  In \cite{Pauls:minimal,DGN,RR,RRas,BC,BSV,DGN:secondvar}, the various authors compute first
variation formulae for $C^2$ smooth noncharacteristic surfaces in the
Heisenberg group.  We note that some of these authors restrict their
attention to certain types of graphs (Euclidean graphs:
\cite{Pauls:minimal,RR}, intrinsic graphs \cite{BSV}), \cite{DGN,DGN:secondvar} deals with
a level set formulation and \cite{BC} provides a completely general
nonparametric first variation formula.  In \cite{CHMY} the authors
compute a first variation formula for $C^2$ noncharacteristic graphs
in any three dimensional pseudo-hermitian space (including, of course,
the Heisenberg group).  In \cite{Selby,Montefalcone} the authors
independently provide a first variation formula for $C^2$
noncharacteristic surfaces in general Carnot groups.  In
\cite{Cole,NS}, the authors provide a first variation
formula for $C^2$ surfaces in Martinet-type space and (2,3) contact
manifolds (respectively).

In \cite{HP2}, we compute the first variation formula for $C^2$
noncharacteristic surfaces in all
so-called {\em vertically rigid} spaces.   To formalize this, we recall some of the basic definitions.
\bgD{SR space}  A {\bf sub-Riemannian (or Carnot-Carath\'eodory) manifold}
is a triple $(M,V_0,\langle \cdot, \cdot \rangle)$ consisting of a
smooth manifold $M$ of dimension $n+1=k+l+1$, a smooth $k+1$-dimensional distribution $V_0
\subset TM$ and a smooth inner product on $V_0$.  This structure is
endowed with a metric structure given by
\[d_{cc} (x,y) = \inf \left \{\int \langle \dot{\gamma},\dot{\gamma}
\rangle^\frac{1}{2} | \gamma(0)=x,\gamma(1)=y, \gamma \in
\mathscr{A}\right \}\]
where $\mathscr{A}$ is the space of all absolutely continuous paths
whose derivatives, when they are defined, lie in $V_0$.
\enD

\bgD{Vertically Rigid} A  {\bf vertical complement} to sub-Riemannian structure is \begin{itemize}
\item a smooth complement $V$ to $V_0$ in $TM$
\item a smooth frame $T_1,\dots T_{l}$ for $V$
\item a Riemannian metric $g$ such that $V$ and $V_0$ are orthogonal, $g$ agrees with $\langle \cdot, \cdot \rangle$ on $V_0$ and $T_1,\dots T_{l}$ are orthonormal.
\end{itemize}
A vertical complement is {\bf rigid} if in addition
\begin{itemize}
\item there exists a partition of $\{1,\dots,l\}$ into equivalence classes such that for all sections $X \in \Gamma(V_0)$, $g ( [X,T_\ua],T_\ub)=0$ if $\ua \sim \ub$.
\end{itemize}

A sub-Riemannian space with a rigid vertically complement is called a
{\bf vertically rigid (VR) space}.

\enD

We note that all of the cases of the first variation formula described above are vertically rigid
and, hence, our first variation formula generalizes all of these.  In
fact (see Theorem 5.5), a slight modification of the argument from
\cite{HP2} allows us to compute a first variation formula in {\em any}
sub-Riemannian manifold:

\begin{MThm}\label{introa1}
Suppose $\Sigma$ is a $C^2$ noncharacteristic hypersurface in a
sub-Riemannian manifold $M$ and
$F$ is a $C^{1;2}$ variation of $\Sigma$ with horizontal variation
function $\vfo $. Then \[\frac{d}{dt}\bigg |_{t=0} P_0(\Sigma_t) = \int_{\Sigma} \vfo  \left(H-\sum\limits_\ub \ap{ [\nu,T_\ub]}{T_\ub}{} \right)\Lambda. 
\]
\end{MThm}
Here, $P_0$ denotes the sub-Riemannian perimeter measure and $H$ denotes the horizontal mean curvature (see Section 3 for a precise definition).  As a
consequence, we have:

\begin{MThm}\label{introa2}
A necessary and sufficient condition for a $C^2$ hypersurface $\Sigma$
to be a noncharacteristic critical point for the horizontal perimeter
measure in the category of $C^1$ hypersurfaces with fixed boundary in
a sub-Riemannian manifold is
\[ \text{div }\nu= H-\sum\limits_\ub \ap{ [\nu,T_\ub]}{T_\ub}{} =0.\]
If the vertical structure is rigid, the second term drops out and the equation becomes
\[ \text{div } \nu=H= 0.\]
\end{MThm}

Recent work of Cheng-Hwang-Yang \cite{CHY} and Ritor\'e-Rosales
\cite{RRas} have shown how to extend the first variation formula to
allow for variations over the characteristic locus in the Heisenberg
group.  Our first main result of this paper is to prove a similar
extension of the previous Theorems to include variations over the
characteristic locus.  

\begin{MThm}\label{introb}
Let $\Sigma$ be a $C^2$ hypersurface in a sub-Riemannian space $M$
with characteristic locus $\Sigma(M)$.  Further suppose that the
Riemannian curvature tensor of $\Sigma$ is bounded and that the
horizontal mean curvature of $\Sigma$, $H$, is in $L^1(\Sigma)$.  
Suppose $F$  is a compactly supported $C^{1;2}$ variation of $\Sigma$ with $F_0$ $C^2$ and variation function $\vf $. Then
\[\begin{split} \frac{d}{dt}\bigg|_{t=0^+} P_0(\Sigma_t) &= \int_{\Sigma \backslash C(\Sigma)} \vf  (\text{div } \nu)  dV_{\Sigma} - \lim_{\ud \to 0} \int_{\partial \Omega_\ud} \vf  \aip{\nu}{N_{\Omega_\ud}}{} dV_{\Omega_\ud} \\
&= \int_{\Sigma \backslash C(\Sigma)} \vf  (\text{div }\nu) - \text{div}_{\Sigma} (\vf  \nu^\top)\;  dV_\Sigma \end{split}\]
Here, $N$ is the Riemannian normal to $\Sigma$, $\nu$ is the unit horizontal normal and $\nu^\top$ is the component of $\nu$ tangential to $\Sigma$.
\end{MThm}

In the previous theorem, the family $\Omega_\ud \subset \Sigma$, $\ud
>0$  satisfies the following conditions:
\begin{itemize}
\item The portion of the boundary $\partial \Omega_\ud$ in the interior of $\Sigma$ is piecewise $C^2$.
\item $C(\Sigma) \subset \Omega_\ud$ for all $\ud>0$
\item $\mu^{n}(\Omega_\ud) \to 0$ as $\ud \to 0$.
\item If we let $N_0$ be the projection of the Riemannian normal to $V$, then as $\ud \to 0$,
 $\int_{\partial \Omega_\ud} |N_0| dV_{\partial \Omega_\ud} \to 0$.
 \end{itemize}
Here $\mu^n$ denotes the $n$-dimensional Riemannian spherical Hausdorff measure.  

We note that due to work of \cite{Balogh,Magnani:cl} on the size of the
characteristic locus, as well as
computation in the appendix, we can show that such a family of sets always exists.  

In keeping with historical terminology, we call critical points of the
perimeter variation {\em minimal surfaces}.

As an application of the general first variation formula, we prove a
version of the Minkowski formula in this setting (see Theorem 6.4 and
Corollary 6.5):

\begin{MThm}\label{introc}
Suppose $\Omega$ is a compact $C^2$ domain with $\Sigma =\partial \Omega$ that is a critical point for perimeter measure with volume constraint. Then
 \[ (Q-1)P_0(\Sigma) =Q H \text{Vol}(\Omega).\]  
\end{MThm}

We note that this formula was shown in groups of Heisenberg type in
\cite{DGN} and in the Heisenberg group in \cite{RRas}.

A number of authors \cite{DGN:secondvar,Selby,Montefalcone} have also computed second variation formulae in the setting of Carnot groups as a tool in the investigation of stable minimal surfaces.  As has been shown recently \cite{DGN:unstable,DGNP,DGNP2,BSV}, stability plays a
crucial role in the study of minimal surfaces in the Heisenberg
group.  Specifically, these papers study analogues of the
sub-Riemannian Bernstein problem and show that without the imposition
of stability on critical points of the first variation of perimeter,
there is no Bernstein-type rigidity.  On the other hand, there is
rigidity in the
presence of the stability condition.  The most general of these
results is an analogue of Riemannian results of
Fischer-Colberie/Schoen \cite{FCS} and Do Carmo/Peng \cite{DCP}:

\begin{Thm}[\cite{DGNP2}]  The only stable $C^2$ complete embeded
  noncharacteristic minimal surfaces without boundary in the first
  Heisenberg group are the vertical planes.
\end{Thm}

To facilitate further study of stable minimal surface, we derive a
second variation formula for $C^2$ noncharacteristic surfaces in
vertically rigid spaces (Theorem 7.4):

\begin{MThm}\label{introd}
Suppose $M$ is a vertically rigid subRiemannian manifold and $F$ is a noncharacteristic $C^{\infty;3}$ variation of $\Sigma\backslash C(\Sigma)$ with compactly supported horizontal variation function. Then
\bgE{Smooth}
\begin{split}
 \frac{d^2}{dt^2}\bigg|_{t=0} P_0( \Sigma_t ) 
 &= \int_{\Xi}   \left[  (\partial_t \vfoc ) H\Lambda \right]_{|t=0} + \int_{\Sigma} \left| \nabla^{0,\Sigma} \vfo  \right|^2 \Lambda  \\
& \quad + \int_{\Sigma} \vfo ^2 \big[ -\text{Ric}^{\nabla}(\nu,\nu) - \text{Tr}(\twoh^2) + H^2\\
 & \qquad   -\ap{\text{Tor}(\nu,e_j)}{T_\ub}{}\ap{\text{Tor}(e_j,T_\ub)}{\nu}{} \\
 & \qquad - \ap{\text{Tor}(e_j,T_\ub)}{N_v}{} \ap{\text{Tor}(\nu,e_j)}{T_\ub}{} \\
 & \qquad  - \Big( e_j \ap{\text{Tor}(\nu,e_j)}{N_v}{} -2\ap{\text{Tor}(\nu,e_j)}{N_v}{} \ap{ \nabla_{e_m} e_j}{e_m}{} \Big)\\
 & \qquad  - \ap{\text{Tor}(\nu,e_j)}{T_\ub}{} (e_j a_\ub) - \ap{\text{Tor}(\nu,e_j)}{N_v}{}^2\Big]\Lambda
 \end{split}\enE
Again, $\nu$ is the unit horizontal normal to $\Sigma$, $\{e_j\}$
form an orthonormal basis for the horizontal tangent space to $\Sigma$,
and the $\{T_i\}$ is the family of vertical vectors fields from the
definition of a vertical complement.  $II_0$ denotes the horizontal second fundamental form. The connection $\nabla$ is adapted to the vertical structure and all torsion terms are associated to $\nabla$.
\end{MThm}

We note that this second variation formula, when restricted to the
special case of the Heisenberg group matches with others in the
literature \cite{DGN:secondvar,DGN,CHY,RRas}.  As a last application, we show that
the conjectural isoperimetric profile in $\mathbb{H}^2$ is indeed a
stable constant mean curvature surface (Section 9).

\section{Notation and conventions}\setS{NC}

To improve economy with the intensive computations throughout this paper, we shall following the following notations and conventions:

\bgEn{\Alph}
 \etem Unless explicitly stated otherwise, roman indices will run from $1 \dots k$, barred roman indices from $0 \dots k$ and greek indices from $1 \dots l$.
 \etem  $\om^{(j)}$ denotes the ordered wedge product of all possible (by index conventions) $1$-forms $\om^i$ with the $j$th form omitted, e.g \begin{align*}
 \om^{(2)} &= \om^1 \wedge \om^3\wedge \dots \wedge \om^k\\ \om^{\bar{2}} &= \om^0 \wedge \om^1 \wedge \om^3  \wedge \dots \wedge \om^k.
 \end{align*}
   We shall also use $\om^{(i,j)}$ with $i<j$ to denote the ordered wedge product with both $i$th and $j$th terms missing and extend to all indices by setting $\om^{(i,j)} =- \om^{(j,i)}$.
\etem (Summation Convention) Whenever the same index appears twice in a term obeying the above conventions, we shall assume that there is an implicit sum over all possible values.
\enEn
\section{SubRiemannian manifolds and vertical structures}\setS{VR} 

We begin with our basic definitions:
\bgD{SR space}  A {\bf sub-Riemannian (or Carnot-Carath\'eodory) manifold}
is a triple $(M,V_0,\langle \cdot, \cdot \rangle)$ consisting of a
smooth manifold $M$ of dimension $n+1=k+l+1$, a smooth $k+1$-dimensional distribution $V_0
\subset TM$ and a smooth inner product on $V_0$.  This structure is
endowed with a metric structure given by
\[d_{cc} (x,y) = \inf \left \{\int \langle \dot{\gamma},\dot{\gamma}
\rangle^\frac{1}{2} | \gamma(0)=x,\gamma(1)=y, \gamma \in
\mathscr{A}\right \}\]
where $\mathscr{A}$ is the space of all absolutely continuous paths
whose derivatives, when they are defined, lie in $V_0$.
\enD

We shall also make the standing assumption that $M$ is oriented.

\bgD{Vertically Rigid} A  {\bf vertical complement} to sub-Riemannian structure is \begin{itemize}
\item a smooth complement $V$ to $V_0$ in $TM$
\item a smooth frame $T_1,\dots T_{l}$ for $V$
\item a Riemannian metric $g$ such that $V$ and $V_0$ are orthogonal, $g$ agrees with $\langle \cdot, \cdot \rangle$ on $V_0$ and $T_1,\dots T_{l}$ are orthonormal.
\end{itemize}
A vertical complement is {\bf rigid} if in addition
\begin{itemize}
\item there exists a partition of $\{1,\dots,l\}$ into equivalence classes such that for all sections $X \in \Gamma(V_0)$, $g ( [X,T_\ua],T_\ub)=0$ if $\ua \sim \ub$.
\end{itemize}

A sub-Riemannian space with a rigid vertically complement is called a
{\bf vertically rigid (VR) space}.

\enD

For convenience of reference, if $M$ is a VR space we decompose
\[ TM = V_0 \oplus \bigoplus\limits_{b \in B}  V_b\]
where $B$ is the set of equivalence classes in $\{1,\dots l\}$ and $V_b = \text{span} \left\{ T_\ub \colon \ub \in b \right\}$. We shall also use $L_b = | b |$. If the vertical structure is not rigid we use the same notation with the understanding that the equivalence relation is simply equality. i.e. $\ua \sim \ub$ if and only if $\ua=\ub$.

Many of the essential computational tools of Riemannian geometry can be generalized or restricted to include subRiemannian geometries with vertical structures.

\bgD{Adapted}
A connection $\nabla$ is adapted to a subRiemannian geometry with vertical structure if
\begin{itemize}
\item $\nabla$ is compatible with $g$
\item $\nabla T_\ub =0$ for $\ub=1\dots l$
\item $\text{Tor}(X,Y)_{|p} \in V_{|p}$ for all sections $X,Y$ of $V_0$
\end{itemize}
\enD

These connections were first defined in \cite{HP2} for VR manifolds but the rigidity assumption is not necessary. The following properties of adapted connections were also proved in \cite{HP2}:

\bgL{Properties} \hfill

\begin{itemize}
\item Every vertical structure admits an adapted connection.
\item For $X,Y$ sections of $V_0$, $\nabla_X Y$ is depends solely on the vertical complement $V$ and not the Riemannian extension $g$.
\item If $\nabla$ is adapted to a rigid vertical structure and $X$ is a section of $V_0$ then
\[ \ap{ \text{Tor}(T_\ua,X) }{T_\ub}{} =0, \qquad \text{if $\ua\sim \ub$.}\]
\end{itemize}
\enL

\bgR{Torsion}
It is the vanishing of these torsion terms that makes VR structures much easier to work with than general subRiemannian manifolds.
\enR

For a $C^1$ hypersurface $\Sigma$ in a subRiemannian manifold, we define the characteristic set of $\Sigma$ to be
\[ C(\Sigma) =\{ p \in \Sigma \colon (V_0)\big |_{p} \subset T_p \Sigma \}.\]

For an oriented $C^1$ hypersurface $\Sigma \subset M$ we define $N$ to be the unit Riemannian normal to $\Sigma$ with respect to $g$ and $N_0$ as the orthogonal projection of $N$ to $V_0$. Away from the characteristic set, we define the {\bf unit horizontal normal} 
\[ \nu = \frac{N_0}{|N_0|}\]

\bgD{HPM}
The horizontal perimeter measure of $\Sigma$ is defined to be
\[ P_0(\Sigma) = \int_\Sigma |N_0| dV_\Sigma \]
where $dV_\Sigma = N \lrcorner dV_g$. 
\enD

For noncharacteristic surfaces, $P_0$ has the alternative descriptions
\[ P_0(\Sigma) = \int_\Sigma \nu \lrcorner dV_\Sigma = \sup\left\{ \int_\Sigma X \lrcorner dV_\Sigma : X \in \Gamma(V_0), |X|=1\right\}\]

There are several natural questions associated to this perimeter measure.

\bgQ{Min}
Among hypersurfaces with the same boundary, which minimizes the horizontal perimeter measure? Can such surfaces be characterized as solutions to a PDE?
\enQ

\bgQ{MinVol}
Among domains of the same volume, which has boundary minimizing horizontal perimeter measure?
\enQ

These problems are studied using variational techniques which may yield critical points rather than true minima. Thus there is another natural question:

\bgQ{Stable}
Of the critical points of horizontal perimeter measure, which are {\bf stable}, i.e. $\frac{d^2}{dt^2}\big |_{t=0} P_0(\Sigma_t) \geq 0$ for any variation of $\Sigma$. 
\enQ

Under the assumptions of rigidity, no characteristic points and $C^2$ regularity, \rfQ{Min} and \rfQ{MinVol} were answered in \cite{HP2} in terms of the horizontal mean curvature.

Suppose $e_0,\dots e_k$ forms a (local) orthonormal frame for $V_0$ such that on $\Sigma \backslash C(\Sigma)$, $e_0$ is the unit horizontal normal to $\Sigma$. Then away from $C(\Sigma)$, the horizontal second fundamental form for $\Sigma$ is defined by
\[ \twoh =  \begin{pmatrix} \ap{ \nabla_{e_1} e_0 }{e_1}{} & \dots &   \ap{ \nabla_{e_1} e_0 }{e_k}{}\\
\vdots & \vdots & \vdots \\
 \ap{ \nabla_{e_k} e_0 }{e_1}{} & \dots &  \ap{ \nabla_{e_k} e_0 }{e_k}{} \end{pmatrix} \]  The horizontal mean curvature is defined by
 \bgE{H}  H = \text{trace} (\twoh).\enE
 We remark that the connection used in these definitions can be either the Levi-Cevita connection for the Riemannian metric or any connection adapted to the vertical structure.
 
 In a VR manifold, $C^2$ minimizers of $P_0$ with  fixed boundary constraint were shown in \cite{HP2} to satisfy $H=0$ away from characteristic points. Likewise $C^2$ minimizers subject to the volume constraint satisfied the condition that the horizontal mean curvature was locally constant away from the characteristic set.

We finish this section by making some remarks on the nature of the equation $H=c$ off the characteristic set. Since any adapted connection is metric compatible, we obtain the following result about the ambient divergence on $M$ from standard results in Riemannian geometry (see \cite{Kobayashi}, appendix 6):
\begin{align*} \text{div }Z &= \text{trace}( \nabla Z + \text{Tor}(Z,\cdot))\end{align*}
Thus if $X_\bt{i}$ denotes a local horizontal orthonormal frame and $Z$ is a horizontal vector field
\begin{align*}
\text{div }Z &= \aip{\nabla_{X_\bt{i}} Z}{X_\bt{i}}{} + \aip{ \text{Tor}(Z,T_\ub)}{T_\ub}{}\\
& = \aip{\nabla_{X_\bt{i}} Z}{X_\bt{i}}{} - \aip{ [Z,T_\ub])}{T_\ub}{}.
\end{align*}
Applying this to $\nu$ yields
\[ \text{div }\nu = H -\aip{ [\nu,T_\ub])}{T_\ub}{}.\]
 In the rigid case, the second term drops out and $H$ naturally takes the form
\[ H = \text{div }\nu .\]  

Elsewhere in the literature, for example \cite{RRas}, the hypersurface divergence has been used instead of the ambient divergence. For completeness we shall now show that the two approaches are equivalent. In particular, our variation formula agrees with that derived  by Rosales and Ritor\'e  for the special case of the first Heisenberg group, . 

\bgD{divE} For vector fields $Z \in \Gamma(TM)$ we define the surface divergence of $Z$  at $p \in \Sigma$  by
\[ ( \divE Z) \omega = \iota^* \mathcal{L}_Z \omega\]
at $p$, 
where $\omega = N \lrcorner dV$ for some extension of $N$ to a unit vector in a neighbourhood of $\Sigma$.
\enD 

We note that it is easy to show that if $Z=Z\upp$ along $\Sigma$ then $\divE Z = \divE Z\upp$.

\bgL{ComputeDiv}
Suppose $\Sigma$ is a $C^2$ hypersurface in a  VR manifold $M$ and $X_1$, \dots $X_n$ is an orthonormal frame for $\Sigma$. Then for any horizontal vector field $Z$
\[ \divE Z = \sum_{j=1}^n \aip{\nabla_{ \iota_* X_j} Z}{\iota_* X_j}{} \]
for any connection $\nabla$ adapted to the VR structure of $M$. 
\enL

\pf We adapt the same argument from  \cite{Kobayashi} as follows: since $\nabla$ is metric compatible we see
\[ \nabla \omega  = \nabla \left( \omega^1 \wedge \dots \wedge \omega^{n-1}\right) = \sum_{j=1}^n \omega^0 \wedge \omega^j_0 \wedge  \omega^{(j)}\]  Thus $\iota^* \nabla \omega = 0$. Therefore if we identify $X_j$ with $\iota_* X_j$,
\begin{align*}
 0&=\left[ \iota^* \left( \nabla_Z - \mathcal{L}_Z \right)  \right] \omega(X_1, \dots X_n)\\
 &= -\iota^* (\mathcal{L}_Z \omega)(X_1,\dots X_n) +   \sum \omega(X_1,\dots (\nabla_Z- \mathcal{L}_Z)X_i, \dots X_n)\\
 &= -(\divE Z) \omega(X_1,\dots, X_n) + \sum_{i=1}^n \omega( X_1 ,\dots \nabla_{X_i} Z +\text{Tor}(Z,X_i), \dots X_n) \\
 &= \left( -\divE Z +  \sum_{j=1}^n \aip{\nabla_{  X_j} Z}{X_j}{}\right) \omega(X_1,\dots X_n)
 \end{align*}
 where the last line follows from the fact that $\aip{\text{Tor}(Z,X_i)}{X_i}{}=0$ whenever $Z$ is horizontal. This last observation fails in the non-rigid general case.
 
\epf

\bgC{divEqual}
Away from the characteristic set $C(\Sigma)$ we have
\[ \text{div } \nu = \divE \nu. \]
\enC

\pf
Away from $C(\Sigma)$ we can choose an orthonormal frame for $\Sigma$ of the form $\{e_j\}_{j=1}^{k-1}$, $\tilde{e}_\ub=  s_\ub \nu +c_\ub T_\ub$ 
where the $e_j$ are horizontal and at least one $c_\ub \ne 0$. Then 
\begin{align*} \divE \nu &= \aip{\nabla_{\tilde{e}_\ub} \nu}{\tilde{e}_\ub}{} + \aip{\nabla_{e_j} \nu}{e_j}{}\\
&= \aip{ \nabla_{\tilde{e}_\ub} \nu}{s_\ub \nu}{} + \text{div } \nu\\
&= \text{div } \nu.
\end{align*}
\epf

\section{Bundles and Variations}\setS{BV}

To describe the variational properties of the horizontal perimeter measure, we shall define a variety of bundles over $M$.

First, we shall denote by $\sM$ the {\bf contact manifold of normalized hypersurface elements}, i.e $\pi\colon \sM \to M$ is the unit tangent bundle over $M$ viewed as a bundle of Riemannian unit normals. We define the $1$-form $\Theta$ on $\sM$ by
\[ \Theta_{|(p,E)}(X) = \ap{\pi_* X}{E}{p} .\]
An immersion $\iota$ of an $n$-dimensional manifold into $\sM$ is said to be \textbf{transverse} if $\pi \circ \iota$ is an immersion and $\iota^* \Theta =0$.

\bgD{N0}
 The function $N_0\colon \sM \to V_0(M)$ is defined by
 \[ N_0(p,E) = \left( (\pi_*E)_0 \right)\big|_{p}  \] 
Here, we use the convention that if $W$ is a vector field on $M$ then $(W)_0$ is its projection to a vector field in $V_0$.

The characteristic slice, $C_\s$ of $\sM$ is the zero level set of $N_0$.
 \enD

There is a natural projection $\pi_\F$ from the Riemannian frame bundle $\F(M)$ to $\sM$ given by
\[ \pi_\F \colon  (p,E_0,\dots E_n) \mapsto (p,E_0).\]
We note that if $E^0, \dots E^n$ denote the tautological forms on $\F(M)$, (i.e. at the point $(p,E_0,\dots E_n)$, $E^j(X) = \ap{\pi_*X}{E_j}{}$, $j=0,\dots,n$) then
\[ \Theta = \sigma^* E^0 \]
for any section $\sigma$ of $\pi_\F$.

\bgD{semibasic}
A differential form $\psi$ on $\sM$ is semibasic if
\[ X \lrcorner \psi =0\]
wherever  $\pi_* X =0$. Thus $\psi$ depends only the the projection to $M$ and the choice of $E_0$.
\enD

For example, it is clear that $\Theta$ is a semibasic $1$-form.

The {\bf bundle of graded orthonormal frames} is the subbundle $\G(M) \subset \F(M)$ such that 
\[ (E_0 \dots E_n) = (e_0\dots e_k, t_{1_1}, \dots, t_{1_{L_1}}, \dots t_{|B|_1},\dots, t_{|B|_{L_{|B|}}}) \] where the $e_j$'s are all horizontal and  $t_{b_1},\dots, t_{b_{L_b}}$ span $V_b$. The reduced structure group of the bundle is then $O(k+1) \times \Pi_{b \in B} O(L_b)$.

For computational purposes it is often easiest to further restrict and insist that vertical frames be exactly those specified in the definition of the vertical structure. If such a vertical structure has been chosen, we shall note this further restriction by $\HF(M)$.

Unfortunately, these graded bundles do not encode enough information to describe the geometry of hypersurfaces of $M$. To compensate for this we also introduce the augmented bundles
\[ \GF[0] = \GF \times \rn{l}, \qquad \HF[,0]= \HF \times \rn{l}.\]
The additional elements will be used to keep track of the dependence of the hypersurface normal directions on the vertical vector fields. 

If an explicit vertical structure has been fixed, there is an alternative presentation of $\sM \backslash C_\s$ that will prove computationally simpler to work with for noncharacteristic variations.

We define the  {\bf contact manifold of horizontally normalized hypersurface elements} to be 
\[ \sMo = \{ (p,Z) \in TM: |Z|=1, Z \in (V_0)_p\} \times \rn{l} .\]
There is a bundle isomorphism $\sM\backslash C_\s \cong \sMo$ given by
\[ (p,E) \mapsto \big(p, |N_0|^{-1} (E)_0, |N_0|^{-1} \ap{E}{T_\ub}{} \big) \]
with inverse
\[ (p,e_0,a_\ub) \mapsto \big(p,  \frac{1}{\sqrt{1+ |a|^2}} (e_0+a_\ub T_\ub) \big).\]
We shall identify $\sMo$ with $\sM \backslash C_\s$ using this bundle isomorphism.

We define a $1$-form $\theta$ on $\sMo$ by
\[ \theta(X) = \ap{ \pi_* X}{e_0 +a_\ub T_\ub}{}\] and note that on $\sMo$ we have
\[ |N_0| = \frac{1}{\sqrt{1+|a|^2}} \text{ and } \Theta = |N_0| \theta.\] 
There is a natural projection from $\pi_\G\colon \GF(M) \to \sMo$ such that if
$ \om^0 \dots \om^k$,  $\eta^1 \dots \eta^{l}$ are the tautological $1$-forms for $\GF(M)$ then 
\[ \theta = \sigma^* (\om^0 +a_\ub \eta^\ub ) \]
for any section $\sigma$ of $\pi_\G$. Since we shall frequently be computing on the frame bundles, we shall often implictly identify $\om^0 +a_\ub \eta^\ub$ with $\theta$.

For the remainder of this section, we shall suppose that $\Sigma$ is an oriented, immersed $C^2$ hypersurface of $M$  realized as the image of the $C^2$ immersion
\[ \iota\colon \Xi \hookrightarrow M\]
for some smooth oriented manifold (possibly with boundary) $\Xi$.

\bgD{Var}
A variation of $\Sigma$ is a map
\[ F\colon \Xi \times (-\e,\e) \to M \] such that
\begin{itemize}
\item Each
$F_t = F(\cdot,t)$ is an immersion of $\Xi$ into $M$.
\item $F_0 = \iota$.
\end{itemize}
The lifted variation $\wh{F} \colon \Xi \times (-\e,\e) \to \sM$ is the map
defined by 
\[\wh{F}(\xi,t)=(F(\xi,t),N\big |_{F(\xi,t)} )\]
where $N$ is the (local) Riemannian unit normal vector to the immersed surface $F_t(\Xi)$ such that the pullback of $N \lrcorner dV$ matches the fixed orientation of $\Xi$.

The variation function of $F$ is $\vf =\vfc (\cdot ,0)$ where $\wh{F}^* \Theta = \vfc dt$. The variation is said to be compactly supported if $\vfc (\cdot,t)$ has compact support for all $t$.
\enD

\bgD{nVar}
When the lifted variation $\wh{F}$ maps into the complement of the characteristic slice $C_\s$, we shall refer to the variation as noncharacteristic. The horizontal variation function for $F$ is then defined by $\vf_0=\vfc _0(\cdot,0)$, $\wh{F}^*\theta = \vfc _0 dt$.
\enD

\bgR{Var}
Since $C_\s$ is closed, if $\Sigma$ has no characteristic points then, shrinking $\e$ if necessary, any variation $F$ will map into $\sMo$. The relationship between the variational functions is just
\[ \vf = |N_0 \circ \wh{F} | \vf_0.\]
\enR

So far, we have not put any regularity conditions on our variations. However, we shall need  precise descriptions of regularity to make our theory optimal.

\bgD{Reg}
The classes of $C^{i;j}$ maps from $\rn[x]{n} \times \rn[t]{}$ to $\rn{}$ for $i,j \geq 0$ are defined inductively by
\begin{itemize}
\item $C^{0;0} =C^0$, i.e. continuous maps.
\item $F \in C^{i+1;0}$ if and only if $F, \pd{F}{x^m} \in C^{i;0}$ for all $m=1 \dots n$. 
\item $F \in C^{0;j+1}$ if and only if $F,\pd{F}{t} \in C^{0;j}$.
\item $F \in C^{i+1;j+1}$ if and only if $F \in C^{i;j+1} \cap C^{i+1;j}$,$\pd{F}{t} \in C^{i+1;j}$  and $\pd{F}{x^m} \in C^{i;j+1}$ for all $m=1 \dots n$.
\end{itemize}
Thus a map is $C^{i;j}$ if up to  $i$ continuous spatial ($x$) derivatives and $j$ temporal ($t$) continuous derivatives can be taken in any order.

Using coordinate charts, this definition extends naturally to define $C^{i;j}$ maps 
\[ \Xi \times \rn{} \to M\]
for smooth manifolds $\Xi$ and $M$.
\enD

\bgR{Incl}
We note in passing that 
\[ C^m = \bigcap\limits_{i+j=m} C^{i;j} \subset C^{m;m}.\]
\enR
The following approximation result will be useful later
\bgL{Approx}
Given a $C^{i;j}$ map $F\colon \Xi \times \rn{} \to M$ that is constant outside $K \times \rn{}$ for some compact set $K \subset \Xi$, there  exists a sequence $F_m$ of $C^{\infty;j}$ maps such that
\bgEn{\arabic}
\etem $F_m$ converges to $F$ in $C^{i;j}$.
\etem If $F(\cdot, t_0) \in C^p$ then $F_m(\cdot,t_0)$ converges to $F(\cdot,t_0)$ in $C^p$.
\enEn 
\enL

This lemma is essentially a version of standard approximation theorems adapted to allow parameters. The reader is referred to \cite{Hirsch} pp.41-55. for a proof that $C^r$ maps between smooth manifolds can be approximated by smooth maps. Theorem 2.3 in \cite{Hirsch} can easily be adapted to give an approximation of $C^{i;j}$ maps from $\Xi \times \rn{} \to \rn{m}$ by $C^{\infty;j}$ maps, with the observation that the mollification process should only be in the $\Xi$ coordinates. Everything else goes through virtually unchanged.

\bgL{vReg}
If $F$ is a $C^{i;j}$ variation $i,j \geq 1$, then $\wh{F}$ is a $C^{i-1,j}$ map.
\enL

\pf
The tangent space to the immersed surface $F_t(\Xi)$  is locally spanned by the vector fields $F_* \partial_{\xi^m}$. Therefore the Gram-Schmidt algorithm followed by a horizontal projection and rescaling, expresses the unit horizontal normal to $F_t(\Xi)$ $\nu$ as a smooth combination of these spanning vector fields.  Thus $\nu$ can be viewed as a $C^{i-1;j}$ function. 

\epf

We now list a few basic regularity properties
\bgL{Sreg}
If $F$ is a $C^{i;j}$ variation and  $\psi$ is a smooth semibasic differential form on $\sM$ then $\wh{F}^* \psi$ is a $C^{i-1;j-1}$ form on $\Xi \times(-\e,\e)$.
\enL

\pf The real issue here is that as a map into $\sM$, $\wh{F}$ only has $C^{i-1;j}$ regularity. However, $\wh{F}^*\psi$ depends tensorially on the projected input $F_*\frac{d}{dt}$, $F_* \frac{d}{d\xi^m}$, which are   $C^{i;j-1}$, $C^{i-1;j}$ vector fields respectively,  and its position $\wh{F}(x,t)$ which is also at least $C^{i-1;j}$.

\epf

In particular, the lemma implies that we can make sense of the pullback of semibasic forms  by $C^{1;2}$ variations despite the fact that the lifted variations are only continuous maps into $\sM$.

\bgC{rReg}
For a $C^{i;j}$ variation $F$ and smooth semibasic form $\psi$, the form \[(\xi,t)\mapsto F_t^* \psi\]
has $C^{i-1;j}$ regularity.
\enC

\pf The proof is identical to the previous lemma except that we no longer need dependence on $F_* \frac{d}{dt}$.

\epf

\bgC{gReg}
For a $C^{i;2}$ variation $F$, $i \geq 1$, the variation function is $C^{i-1}$.
\enC

Furthermore, it will be of interest to note that, locally at least, every function on $\Sigma$ can be realized as a variation function.

\bgL{realize}
For every point $p =\iota(\xi) \in \Sigma$ there exists a neighborhood $\xi \in U \subset \Xi$ such that every $C^j$ function $\vf$, $j=1,2$, on $U$ is the restriction of the variation function for a $C^{j;\infty}$ variation of $\Sigma$.
\enL

This is shown using standard arguments with Pfaff coordinates (see \cite{BGG1}, p.16). The restriction $j=1,2$ is due to the fact that $\Sigma$ is only assumed to be $C^2$.

\section{Horizontal Perimeter Measure and the First Variation}\setS{FV}

Throughout this section we suppose $M$ is a subRiemannian manifold with a chosen vertical structure.

Recall that on $\pi\colon \GF(M) \to M$ we have the tautological $1$-forms $\om^j$ and $\eta^j$ by
\begin{align*}
\om^j (X) &= \ap{ \pi_*X}{e_j}{}, \qquad j =0 \dots k \\
\eta^j (X) &= \ap{ \pi_*X}{t_j}{}, \qquad j =1 \dots n-k \\
\end{align*}
When computing locally with a noncharacteristic variation,  we shall use the pullbacks of the forms $(\sigma \circ \wh{F})^* \om^1$,$\dots$, $(\sigma \circ \wh{F})^* \om^k$, $(\sigma \circ \wh{F})^*\eta^1 \dots$ $(\sigma \circ \wh{F})^*\eta^{n-k}$ together with $dt$ as a frame on $\Xi \times (-e,\e)$. Here $\sigma$ is any section of $\pi_\G$. We shall drop the $(\sigma \circ \wh{F})^*$'s when referring to this framing.

\bgR{Lambda}
We shall refer to the pullback of $\Lambda$ to $\sMo$ also as $\Lambda$ and to the pullback to $\Xi$ of $\Lambda$ by $\wh{F}_t$ as  $\Lambda_t$. \enR

\bgL{AP}
If $\Sigma$ is a noncharacteristic hypersurface of $M$ and $f$ is a tranvserse immersion $\Sigma \hookrightarrow \sMo$  then
\[ P_0(\Sigma) = \int_\Sigma f^*\Lambda \]
\enL

\pf  This is obvious from the definitions. 

\epf

Most computations will be undertaken on the frame bundles. In addition to the tautological forms, we also have the bundle structural equations for the connection (see \cite{HP2})
\bgE{Structure}
\begin{split}
d\om^\mb{j} &= \om^\mb{m} \wedge \om_\mb{m}^\mb{j} + \tau^\mb{j}\\
d\eta^\ub &= \eta^\ua \wedge \eta^\ub_\ua + \tilde{\tau}^\ub.
\end{split}
\enE 
When working on $\HF[,0]$ (as opposed to $\GF[0]$) the vertical torsion pieces take the form
\[ \tilde{\tau}^\ub = d \eta^\ub.\]
As was noted in \cite{HP2}, if the vertical structure is rigid then $\tilde{\tau}^\ub(X,t_\ub)=0$.

\bgD{Psi}
We define the $1$-form $\Psi$ on $\GF[0]$ by
\[ \Psi = (-1)^{j-1} \om^j_0 \wedge  \omj{j} \wedge \eta^*  \]
\enD

\bgL{Old}
If $\sigma$ is any section of $\pi_\G$ then 
\[ d\Lambda = \theta \wedge \sigma^* \Psi + \sigma^* \left( (-1)^{k+\ub-1} \om^* \wedge \tilde{\tau}^\ub \wedge \etj{\ub} \right). \] 
\enL

\pf
This was essentially proved in \cite{HP2} section 4, but without the assumption of rigidity the torsion term cannot be dropped.

\epf

The importance of $\Psi$ lies in \rfL{Old} and the fact that for a noncharacteristic variation of $\Sigma$,
\[ (\sigma \circ \wh{F}_0)^* \Psi = H \Lambda_0. \]

\bgT{FirstVariation}
Suppose $\Sigma$ is a $C^2$ noncharacteristic hypersurface in $M$ and $F$ is a $C^{1;2}$ variation of $\Sigma$ with horizontal variation function $\vfo $. Then \[\frac{d}{dt}\bigg |_{t=0} P_0(\Sigma_t) = \int_{\Sigma} \vfo  \left(H-\sum\limits_\ub \ap{ [\nu,T_\ub]}{T_\ub}{} \right)\Lambda. 
\]

\enT

\pf
For the rigid case, this was proved in \cite{HP2}, section 4 . A very minor modification using \rfL{Old} without dropping the torsion term gives the general case. For completeness, we sketch out the argument: first note
\[\frac{d}{dt}\bigg |_{t=0} P_0(\Sigma_t) = \left[ \int_{\Xi} \Lie{\partial_t} \wh{F}^* \Lambda \right]\bigg |_{t=0} \]
Now $\Lambda$ is a spatial form on $\sMo$, so we can only guarantee $C^{0;1}$ of $\wh{F}^* \Lambda$ on $\Xi \times (-\e,\e)$. However, since the variation has $C^{1;2}$  regularity and $\Sigma$ itself is a $C^2$ hypersurface, at $t=0$ we can differentiate on  $\HF[,0]$ to see
\[ \begin{split} \frac{d}{dt}\bigg |_{t=0} P_0(\Sigma_t) &=  \left[ \int_{\Xi} d( \partial_t \lrcorner \wh{F}^* \Lambda)\right]\bigg |_{t=0} +\left[ \int_{\Xi} \partial_t \lrcorner \wh{F}^* d\Lambda \right]\bigg |_{t=0}\\
&= \int_{\partial \Xi} \left[ \partial_t  \lrcorner \wh{F}^* \Lambda  \right]\bigg |_{t=0} \\ & \qquad + \int_{\Xi} \partial_t \lrcorner (\sigma \circ \wh{F})^* \left(\theta \wedge \Psi +\left( (-1)^{k+\ub-1} \om^* \wedge \tilde{\tau}^\ub \wedge \etj{\ub} \right) \right)\bigg |_{t=0}\\
&=\int_{\Xi} \vfo \left[  (\sigma \circ \wh{F}_0)^* \Psi +  (\sigma \circ \wh{F})^* \left(- \aip{[\nu,T_\ub]}{T_\ub}{} \om^* \wedge \eta^* \right) \right]\\
&= \int_{\Sigma} \vfo  \left(H-\sum\limits_\ub \ap{ [\nu,T_\ub]}{T_\ub}{} \right)\Lambda. \end{split}\]
Where we have used the fact that
\[ \tilde{\tau}^\ub (\nu,T_\ub) = -\aip{[\nu,T_\ub]}{T_\ub}{}.\]
\epf

\bgC{Perimeter Critical}
A necessary and sufficient condition for a $C^2$ hypersurface $\Sigma$ to be a noncharacteristic critical point for the horizontal perimeter measure in the category of $C^1$ hypersurfaces with fixed boundary is
\[ \text{div }\nu= H-\sum\limits_\ub \ap{ [\nu,T_\ub]}{T_\ub}{} =0.\]
If the vertical structure is rigid, the second term drops out and the equation becomes
\[ \text{div } \nu=H= 0.\]
\enC

\bgR{Citations}
This is the first result of this nature for completely general
subRiemannian manifolds. The rigid case was shown in \cite{HP2}.
Prior results include numerous important cases:  level sets in Carnot groups
\cite{DGN,Montefalcone}, for graphs in three dimensional strictly
pseudoconvex pseudohermitian manifolds \cite{CHMY}, Martinet-type
spaces \cite{Cole}, for graphs in the Heisenberg group
\cite{Pauls:minimal,GP}, for parametrized surfaces in the Heisenberg
group \cite{BC}, for intrinsic graphs in the Heisenberg groups
\cite{BSV}, and for surfaces in (2,3) contact manifolds \cite{NS}.
\enR

 In the presence of characteristic points, the situation becomes more complicated. It is to this case that we now direct our attention. To avoid needless repetition, we shall make the following assumptions throughout this section.
\bgEn[Conditions]{\Alph}
\etem $\Sigma$ is an oriented $C^2$ hypersurface with piecewise $C^1$ boundary  in some $n+1$-dimensional VR manifold $M$.
\etem The Riemannian unit normal to $\Sigma$ will be denoted $N$. Off
the characteristic set $C(\Sigma)$, the unit horizontal normal $\nu =
\frac{1}{|N_0|} N_0$.
\etem $\Sigma$ is the image of the immersion $\iota \colon \Xi \to M$ with $\Xi \subset \rn{n}$.
\etem $F \colon \Xi \times (-\e,\e) \to M$ is a $C^{1;2}$ variation of $\Sigma$ with $F_0$ a $C^2$ mapping. In particular this implies that $\wh{F}$ is a $C^{0;2}$ map.
\etem The horizontal mean curvature of $\Sigma$, $H \in L^1(\Sigma)$.
\etem The Riemannian curvature tensor of $\Sigma$ is bounded.
\enEn

\bgR{} For any transverse immersion $f$  of $\Sigma$ into $\sMo$, the $N_0$ referred to above is equivalent to $f^* N_0$. Likewise we can pull $N_0$ back to $\Xi \times (-\e,\e)$ and we shall not make any notational distinction between them. 
\enR

The necessary observation for studying variations for hypersurfaces with characteristic points is the following: 

Suppose $F$ is a variation of $\Sigma$. Set $C(\Xi) = \iota^{-1}C(\Sigma)$ and note that $C(\Xi)$ is a closed subset of $\Xi$. Furthermore by the results of the appendix, $C(\Xi)$ must have Hausdorff dimension $\leq n-1$. Let $U$ be any open subset of $\Xi$ containing $C(\Xi)$. By shrinking $\e$ if necessary, $F$ induces a noncharacteristic variation $F_H$  of $\Sigma \backslash \iota(U)$ as discussed in \rfS{BV}.  Furthermore, if $\vf $ is the variation function for $F$, then $|N_0|^{-1}\vf $ is the variation function for $F_H$. In particular
\[ (\partial_t \lrcorner \wh{F}^* \Theta )_{|t=0} dV_{\Sigma} = (\partial_t \lrcorner \wh{F}_H^* \theta)_{|t=0} \Lambda.\]

Before diving into the general first variation formula, we shall need some technical lemmas.

\bgL{N0}
With the assumptions listed above,
\begin{itemize}
\item $|N_0(\cdot, 0)|$ is a Lipschitz function on $\Xi \times \{t\}$.
\item $|N_0(\cdot, t)|$ is $C^{0;1}$ off $C(\Sigma_t)$ and has bounded distributional temporal derivative on all of $\Sigma$.
\item The one-sided derivative $\frac{d}{dt}\bigg |_{t=0^+} |N_0(t)|$ exists everywhere, is continuous off $C(\Sigma)$ and is bounded on $\Sigma$ 
 \end{itemize}
\enL

\pf The first part follows from the fact that $N_0(\cdot,0)$ is a $C^1$ map as $F_0$ is $C^2$. The  remaining parts of the lemma are obvious properties of the absolute value of a $C^1$ function from $\rn{}$ to $\rn{}$.

\epf

\bgL{Curve}
Suppose $U$ is an open set in $\Sigma$ such that $\partial U$ does not intersect $C(\Sigma)$. Then for any vector field $X$ on $M$,
\[  X \lrcorner \Lambda _{|\partial U} =|N_0| \aip{X}{N_{\partial U}}{} dV_{\partial U}- \aip{\nu}{N_{\partial U}}{} \aip{X}{N}{} dV_{\partial U}\]
where $N_{\partial U}$ is the Riemannian unit normal inside $\Sigma$ to $\partial U$.
\enL

\pf
Note that off $C(\Sigma)$, $N = |N_0| \nu +|N_0| a_\ub T_\ub $ for constants $a_\ub$. The last piece can be rewritten as $\tilde{a} \tilde{T}$ for some unit vector $\tilde{T}$ orthogonal to $V_0$. Thus we can construct a vector field $\tilde{e}$ along $\Sigma\backslash C(\Sigma)$ by $\tilde{e} = \tilde{a} \nu - |N_0| \tilde{T}$. Then since we must have $|N_0|^2 + \tilde{a}^2 =1$, clearly $\nu = |N_0| N + \tilde{a} \tilde{e}$.   

Now $X \lrcorner \Lambda = X \lrcorner \nu \lrcorner dV$, thus splitting $\nu$ into pieces orthogonal and tangent to $\Sigma$ we have
\[ X \lrcorner \Lambda = |N_0| X \lrcorner dV_\Sigma + \beta X \lrcorner  \tilde{e} \lrcorner dV\]
Pulling back to $\partial U$ immediately yields
\[ X \lrcorner \Lambda _{|\partial U} = |N_0| \aip{X}{N_{\partial U}}{} dV_{\partial U} - \tilde{a} \aip{X}{N}{}\aip{\tilde{e}}{N_{\partial U}}{} dV_{\partial U}.\]
Noting that $\tilde{a} \tilde{e}$ is the tangential component of $\nu$ then completes the proof.

\epf

\bgL{Technical}
There exists a family $\Omega_\ud \subset \Sigma$, $\ud >0$ 
such that
\begin{itemize}
\item The portion of the boundary $\partial \Omega_\ud$ in the interior of $\Sigma$ is piecewise $C^2$.
\item $C(\Sigma) \subset \Omega_\ud$ for all $\ud>0$
\item $\mu^{n}(\Omega_\ud) \to 0$ as $\ud \to 0$.
\item $\int_{\partial \Omega_\ud} |N_0| dV_{\partial \Omega_\ud} \to 0$ as $\ud \to 0$.
\end{itemize}
where $\mu^n$ is the $n$-dimensional Riemannian spherical Hausdorff measure.
\enL

\pf
By \rfT[CS]{Main}, the characteristic set $C(\Sigma)$ is compact and has Hausdorff dimension $\leq n-1$. Thus for any $\ud, p>0$ we can construct a finite collection of Riemannian balls $\mathcal{U}$ of radius $\e$ covering $C(\Sigma)$  such that 
\[ \sum\limits_{\mathcal{U}}  \e^{n-1+p}  < \ud. \]
For each $\ud>0$ take $\Omega_\ud = \Sigma \cap \bigcup\limits_{\mathcal{U}} B$ for some $0<p <<1$. This family clearly satisfies the first three properties. 

The standing assumption on the Riemannian curvature tensor (F) implies that  for some constant $C$
\[ \int_{\partial \Omega_\ud} |N_0| dV_{\Omega_\ud} \leq \e \int_{\partial \Sigma} dV_{\Omega_\ud} + C \sum \limits_{\mathcal{U}} \e^{n-1+1} \]
which tends to zero as $\ud \to 0$.

\epf

We are now in a position to state and proof the main result of this section, the first variation formula for perimeter measure of $C^2$ surfaces.


\bgT{Main}
Suppose $F$  is a compactly supported $C^{1;2}$ variation of $\Sigma$ with $F_0$ $C^2$ and variation function $\vf $. Then
\[\begin{split} \frac{d}{dt}\bigg |_{ t=0^+} P_0(\Sigma_t) &= \int_{\Sigma \backslash C(\Sigma)} \vf  (\text{div } \nu)  dV_{\Sigma} - \lim_{\ud \to 0} \int_{\partial \Omega_\ud} \vf  \aip{\nu}{N_{\Omega_\ud}}{} dV_{\Omega_\ud} \\
&= \int_{\Sigma \backslash C(\Sigma)} \vf  (\text{div }\nu) - \text{div}_{\Sigma} (\vf  \nu^\top)\;  dV_\Sigma \end{split}\]
where $\Omega_\ud$ is any family satisfying the conditions of \rfL{Technical} and $\nu^\top$ is the Riemmanian orthogonal projection of $\nu$ onto $T\Sigma$.
\enT

\pf Using the existence of a family of neighborhoods of $C(\Sigma)$ as in \rfL{Technical} and pulling back to $\Xi$, we can immediately decompose
\bgE{Main}
\frac{d}{dt}\bigg |_{ t=0^+} P_0(\Sigma_t) = \left[ \int_{\Xi
    \backslash \Xi_\ud } \mathcal{L}_{\partial_t}
  (\wh{F}^*\Lambda)\right]\bigg |_{t=0}+ \frac{d}{dt}\bigg |_{t=0^+}  \int_{\Omega_\ud}  \left( |N_0(t)| dV_{\Sigma_t} \right).
\enE
Using the results of \cite{HP2} as discussed in \rfT{FirstVariation} we can reduce the first term to 
\[  \int_{\Sigma \backslash \Omega_\ud} \vf  (\text{div }\nu)
dV_\Sigma  + \left[ \int_{\partial \Xi_\ud} \partial_t \lrcorner
  \wh{F}^*\Lambda \right]\bigg |_{t=0}. \]
This equals 
\[\begin{split}
 \int_{\Sigma \backslash \Omega_\ud} &\vf  (\text{div }\nu)  dV_\Sigma
 +\int_{\partial \Omega_\ud} |N_0| \aip{\wh{F}_* \frac{d}{dt}\bigg |_{t=0}}{N_{\partial \Omega_\ud}}{} dV_{\partial \Omega_\ud} \\
 & \qquad - \int_{\partial \Omega_\ud} \vf  \aip{\nu^\top}{N_{\partial \Omega_\ud}}{} dV_{\partial \Omega_\ud}
 \end{split}
 \]
by \rfL{Curve}. However the middle term can be neglected as $|N_0|$ is Lipschitz on $\Sigma$, with the other terms bounded, and so the integral will vanish as $\ud \to 0$ by \rfL{Technical}. Thus we need only consider the contribution of 
\[ \int_{\Sigma \backslash \Omega_\ud} \vf  (\text{div }\nu)  dV_\Sigma - \int_{\partial \Omega_\ud} \vf  \aip{\nu^\top}{N_{\partial \Omega_\ud}}{} dV_{\partial \Omega_\ud}\] which by the Riemannian divergence theorem can also be expressed as
\[ \int_{\Sigma \backslash \Omega_\ud} \vf  (\text{div }\nu)   - \text{div}_\Sigma (\vf  \nu^\top) \; dV_\Sigma.\]
Now the second term of \rfE{Main} decomposes as
\[ \int_{\Omega_\ud}  \left(  \frac{d}{dt}\bigg |_{t=0^+} |N_0(t)| \right)dV_{\Sigma} + \int_{\Xi_\ud} |N_0| \mathcal{L}_{\partial_t} \iota^*_t dV_{\Sigma_t}\]
By \rfL{N0},  $|N_0|$ is Lipschitz on $\Sigma$ and vanishes on
$C(\Sigma)$. The second of  integral can therefore be uniformly
bounded by a fixed constant times $\ud\int_{\Xi_\ud} \mathcal{L}_{\partial_t} \iota^*_t dV_{\Sigma_t}$. Furthermore  by \rfL{N0} again, $|N_0|$ has bounded distributional derivative, so the first integral is bounded by $\mu^{n}(\Omega_\ud)$. As $\wh{F}$ is $C^1$, $\int_{\Xi_\ud} \mathcal{L}_{\partial_t} \iota^*_t dV_{\Sigma_t} $ is bounded on $\Xi$. Therefore as $\ud \to 0$ the second term of \rfE{Main} tends to zero. Therefore letting $\ud \to 0$ yields the desired result.

\epf

\bgC{Minimizer}
A necessary and sufficient condition for a $C^2$ surface to be a critical point of  horizontal perimeter measure in the category of $C^2$ hypersurfaces with  fixed boundary is \[\text{div } \nu = H- \sum\limits_\ub \ap{ [\nu,T_\ub]}{T_\ub}{}  =0\] on $\Sigma \backslash C(\Sigma)$ and \bgE{Min}  \lim_{\ud \to 0} \int_{\partial \Omega_\ud} \vf  \aip{\nu}{N_{\Omega_\ud}}{} dV_{\Omega_\ud} =- \int_{\Sigma \backslash C(\Sigma)}  \text{div}_{\Sigma} (\vf  \nu^\top) dV_\Sigma=0\enE
for all  compactly supported $C^1$ functions $\vf $.
\enC

\pf This follows from the fact that every compactly supported $C^2$
function can be realized as the variation function of a $C^2$
variation of $\Sigma$ and our assumption (E) that the horizontal mean
curvature is in $L^1$.
\epf
 
\bgC{Isoperimetric}
A $C^2$ perimeter critical domain in the category of $C^2$ domains with volume constraint must have boundary $\Sigma$ satisfying both $\text{div }\nu =c$ off  $C(\Sigma)$ for some constant $c$ and \rfE{Min}.
\enC

\pf
Take any open $\tilde{U} \subset M$ small enough so that $dV=d\mu$ is exact on $\tilde{U}$ and then set $U = \tilde{U} \cap \Sigma$. Then any variation with variation function $\vf $ supported inside $U$ must satisfy
\[ \frac{d}{dt}\bigg |_{ t=0} \left( P_0(\Sigma_t) - \int_{\Sigma_t} c \mu  \right) =0 \]
for some constant $c$. But  
\begin{align*}
\frac{d}{dt}\bigg |_{t=0} \int_{\Sigma_t} c \mu  &= \int_{\Sigma} c \Lie{\partial_t} \mu = \int_{\Sigma} \partial_t \lrcorner dV + d(\partial_t \lrcorner \mu)\\ 
&= \int_{\Sigma} \partial_t \lrcorner \Theta \wedge dV_\Sigma= \int_\Sigma c\vf    dV_\Sigma.
\end{align*}
Thus by using a partition of unity we have that
\[ \int_{\Sigma \backslash C(\Sigma)} \left( \text{div }\nu-c\right)  dV_\Sigma - \lim_{\ud \to 0} \int_{\partial \Omega_\ud} \vf  \aip{\nu}{N}{} dV_{\partial \Omega_\ud} \]
must vanish for all $C^2$ functions $\vf $ on $\Sigma$. Here we can use the same constant $c$ on each supporting patch of the partition as the constants must agree on overlaps.

\epf

\bgR{Char}
For any $C^2$ hypersurface $\Sigma$ such that the characteristic set has Hausdorff dimension $<n-1$, the family $\Omega_\ud$ can be chosen so that condition \rfE{Min} is automatically satisfied. This follows easily from the observation that if we follow the construction of \rfL{Technical} then $\int_{\partial \Omega_\ud} dV_{\partial \Omega_\ud} \to 0$.  For example, as seen in \cite{Balogh,Magnani,CHY}, in the  Heisenberg groups $\hn{m}$ the characteristic set of any $C^2$ hypersurface has dimension $\leq m$. Thus for $m >1$ there is no constraint on the characteristic set of minimal surfaces.
\enR

\bgR{nuCts}
Let $\Sigma$ be a critical point for perimeter variation (with or without volume constraint). Suppose that $p\in C(\Sigma)$ and in a small neighborhood $U$ of $p$, $C(\Sigma)$ is an embedded submanifold of dimension $n-1$ dividing $U$ into two regions $U^+$ and $U^-$. If we further suppose that $\nu$ extends continuously to $\nu^+$, $\nu^-$ on the boundaries of $U^+$ and $U^-$ respectively, then 
\[ \aip{\nu^+}{N_c}{} - \aip{\nu^-}{N_c}{}=0\]
where $N_c$ is the normal to $C(\Sigma)$ in $U$ pointing into $U^-$.

This follows immediately from the divergence integral form of
\rfE{Min}.  Note that, when restricted to the Heisenberg group, this
was observed in \cite{CHMY,RRas} where the structure of the
characteristic locus is known to be either lower dimensional or of
this form.
\enR

 \section{Application: a Minkowski formula for CMC surfaces}\setS{MF}
  
Throughout this section we shall suppose $M$ has a globally defined rigid vertical complement with decompostion
\[ TM = V_0 \oplus \bigoplus\limits_{b \in B} V_b\]
as in \rfS{VR}. 

\bgD{Dilation} A dilating flow for a subRiemannian manifold with global vertical complement is a map $D\colon M \times \rn{} \longrightarrow M$ and constants $\gamma_\ub$, $\ub =1..l$ such that
 \begin{itemize}
  \item $(D_\lambda)_*$ maps $V_0$ to $V_0$.
  \item $\aip{ (D_\lambda)_* Y}{(D_\lambda)_* Z}{D_\lambda p} = e^{2\lambda} \aip{Y}{Z}{p}$ for all sections $Y,Z$ of $V_0$. 
  \item $(D_\lambda)_* T_\ub =e^{\gamma_\ub \lambda} T_\ub$ for all $\ub$. 
 \end{itemize}
 Associated to a dilating flow are the dilation operators defined by
 \[ \delta_\lambda = D_{\log{\lambda}}\]
 and the generating vector field $X$ defined by
 \[ X_p = \frac{d}{d\lambda}_{|\lambda=0} D_\lambda(p).\]
 The homogeneous dimension of $M$ is given by \[ Q = k +1+ \sum_{\ub=1}^l \gamma_\ub.\]
 \enD
 
 For compactness of notation, we shall write $\lambda p$ for $D_\lambda (p)$ and $\lambda_* Y$ for $(D_\lambda)_* Y$.
 
 A dilating flow naturally lifts to a global flow $\wh{D}$ on the contact bundle $\sMo$.  If  $\wh{p} = (p,\nu,a_\ub)$ with $\ub=1..l$ then 
 \bgE{lift} \lambda \wh{p}= \left( \lambda p, e^{-\lambda} \lambda_* \nu, e^{(1-\gamma_\ub )\lambda} a_\ub  \right).\enE
This lifts ensures that the middle term is still unit length and that if $\nu +a_\ub T_\ub $ is a normal vector for the surface $\Sigma= \{\phi =0\}$, then $e^{-\lambda} \lambda_* \nu+ e^{(1-\gamma_\ub )\lambda} a_\ub T_\ub $ is a normal vector to $\Sigma_\lambda= \{D_{-\lambda}^* \phi =0\}$. We also note that $\pi \circ \wh{D}_\ul = D_\ul \circ \pi$. The generator of the lifted flow will be denoted \wh{X}.

\bgL{Invariant}
The contact form $\theta_\wh{p} (\wh{Y}) = \aip{ \pi_* \wh{Y}}{\nu+a_\ub T_\ub }{p}$ has the property
\[ \mathcal{L}_\wh{X} \theta = \theta.\]  
\enL

\pf
We compute
 \begin{align*}
 ( \wh{D}_\lambda^* \theta)_\wh{p}(\wh{Y}) &= \theta_{\lambda \wh{p}} (\lambda_* \wh{Y}) =  \aip{ \pi_* \lambda_* \wh{Y}}{ e^{-\lambda} \lambda_* \nu + e^{(1-\gamma_\ub  )\lambda} a_\ub  T_\ub }{\lambda p}\\
 &= e^{\lambda} \aip{ \pi_* \wh{Y}}{ \nu+a_\ub T_\ub }{p}= e^{\lambda} \theta_\wh{p}(\wh{Y}).
\end{align*}
The result is a direct consequence.

\epf

\bgL{LieL}
The horizontal perimeter measure form $\Lambda$ on $\sMo$ has the following dilation property:
\[ \Lie{\wh{X}} \Lambda = (Q-1) \Lambda.\]
\enL

\pf We first note that clearly
\[ D_\lambda^* dV =  e^{\lambda Q} dV.\]
Now
\begin{align*}
( \wh{D}_\lambda^* \Lambda  )_\wh{p} \left(\wh{Y}_1, \dots, \wh{Y}_n\right) &= \Lambda_{ \lambda \wh{p}} \left(\lambda_* \wh{Y}_1, \dots \lambda_* \wh{Y}_n \right)\\
&=\pi^* \left( e^{-\lambda} \lambda_* \nu \lrcorner dV \right)_{\lambda \wh{p} }  \left( \lambda_* \wh{Y}_1, \dots \lambda_*  \wh{Y}_n \right)\\
&= e^{-\lambda} dV_{\lambda p}  \left( \lambda_* \nu, \lambda_* \pi_* \wh{Y}_1, \dots, \lambda_* \pi_* \wh{Y}_n \right)\\
&= e^{-\lambda} (D_\lambda)^* dV_p \left(\nu, \pi_* \wh{Y}_1,\dots \pi_* \wh{Y}_n \right)\\
&= e^{(Q-1)\lambda} dV_{ p} \left(\nu, \pi_* \wh{Y}_1,\dots \pi_* \wh{Y}_n \right)\\
&= e^{(Q-1)\lambda} \Lambda_{\wh{p}}\left( \wh{Y}_1, \dots \wh{Y}_n \right)
\end{align*}
The result immediately follows

\epf

In the presence of a dilation, we can define 
\[ \Upsilon = Q^{-1} X \lrcorner dV \]
so that $d\Upsilon =Q^{-1} \Lie{X} dV = dV$
 and $\Lie{X} \Upsilon = Q\Upsilon$. We also define $\wh{\Upsilon} =\pi^* \Upsilon$, the pullback of $\Upsilon$ to $\sMo$. Since $\wh{D}_\ul^* \wh{\Upsilon} = \pi^* D_\ul^* \Upsilon$, we immediately see that $\Lie{\wh{X}} \wh{\Upsilon} = Q \wh{\Upsilon}$.

Now suppose $\Sigma$  is a noncharacteristic $C^2$ hypersurface of $M$ with constant mean curvature $H$. Then $\Sigma$ embeds naturally  as $\wh{\Sigma}$ into $\sMo$ and
\[ \int_\wh{\Sigma}  \Lie{\wh{X}} ( \Lambda - H \wh{\Upsilon} )  = \int_\wh{\Sigma} (Q-1)\Lambda-QH \wh{\Upsilon}. \]
But since $\pi^* dV = \theta \wedge \Lambda$ and $d\Lambda = H\pi^* dV$
\bgE{Minkowski}\begin{split} 
 \int_\wh{\Sigma} \Lie{\wh{X}} ( \Lambda - H\wh{\Upsilon} ) &= \int_{\partial \wh{\Sigma}}  \wh{X} \lrcorner \Lambda  + \int_{\wh{\Sigma}} \wh{X} \lrcorner (H\pi^*dV)  -HQ\wh{\Upsilon} \\
& =\int_{\partial \wh{\Sigma}} \wh{X}  \lrcorner\Lambda.
\end{split}
\enE 
After pulling back to $\Sigma$ along the natural inclusion into $\sMo$, we have now established a Minkowski type identity for $C^2$ noncharacteristic patches. Namely
\bgE{NonC} (Q-1) P_0 (\Sigma) = (Q-1) \int_\Sigma \Lambda= QH \int_\Sigma \Upsilon + \int_{\partial \Sigma} X \lrcorner \Lambda \enE

 \bgT{Main}
 Suppose $\Sigma$ is a  $C^2$ hypersurface with piecewise $C^1$ boundary such that $H$ is constant off $C(\Sigma)$ and $\Sigma$ satisfies the constraint \rfE[FV]{Min}. Then
 \[ (Q-1) P_0(\Sigma) = Q \int_\Sigma H \Upsilon + \int_{\partial \Sigma} X \lrcorner \Lambda.\] 
 
 \enT
 
 \pf
As before we use the family of open sets $\Omega_\ud$ containing $C(\Sigma)$ constructed in \rfL[FV]{Technical} and set $\Sigma_\ud = \Sigma - \Omega_\ud$. 
 
 Then by \rfE{NonC}  we know that
 \bgE{ud} (Q-1) P_0(\Sigma_\ud) = Q \int_{\Sigma_\ud} H \Upsilon + \int_{\partial \Sigma_\ud} X \lrcorner \Lambda.\enE
 But by an argument identical to \rfT[FV]{Main}, the internal portions of boundary integral will tend to zero as $\ud \to 0$, leaving the desired equality.  
  
 \epf
  
 \bgC{FV}
 Suppose $\Omega$ is a compact $C^2$ domain with $\Sigma =\partial \Omega$ that is a critical point for perimeter measure with volume constraint. Then
 \[ (Q-1)P_0(\Sigma) =Q H \text{Vol}(\Omega).\]  
 \enC
 
 For the Heisenberg groups, this result was first shown in \cite{RRas}.

\section{The Second Variation}\setS{SV}

We shall now attempt the arduous task of describing a general second variation formula under the assumption of rigidity. This unfortunately is just a long tedious exercise in computing derivatives using the structural equations of the adapted connection in $\HF[,0](M)$. The underlying idea is differentiate on both the frame bundle and on $\Xi \times (-\e,\e)$ and compare results.

To aid with the long computations to follow, we shall briefly list the standing assumptions and notational conventions of this section. To save time and space, we shall also adopt the habit of absorbing all unnecessary terms that do not affect the relevant computations into ``junk'' collections.

\bgEn{\Alph}
 \etem $M$ is a vertically rigid subRiemannian manifold of dimension $n+1=k+1+l$. 
  \etem Unless otherwise specified $\Sigma= \iota(\Xi)$ is a $C^{\infty}$ hypersurface in $M$ with no characteristic points.
 \etem $F\colon \Xi \times (-\e,\e)$ is a $C^{\infty;3}$ noncharacteristic variation with $\wh{F}^*\theta =\vfoc $ and $\vfo =\vfoc (\cdot,0)$.
 \etem And recall: roman indices run from $1 \dots k$, barred roman indices from $0 \dots k$ and greek indices from $1 \dots l$.
\enEn

We use the framing $dt, \wh{F}^* \om^j, \wh{F}^*\eta^\ua$ on $\Xi \times (-\e,\e)$ and will generally omit the $\wh{F}^*$. We shall use the notation $\om^* = \om^1 \wedge \dots \wedge \om^k$ and $\eta^* = \eta^1 \wedge \dots \wedge \eta^l$.

We define a variety of tensors by pulling back the structural equations to $\Xi \times (-\e,\e)$.
\[ \begin{split} 
\wh{F}^* \om^\mb{i}_\mb{j}  &= \uc^\mb{i}_\mb{j} dt + \Gm^\mb{i}_{\mb{j} m} \om^m +  \Gamma^{\mb{i}}_{\mb{j}\ua}  \eta^\ua\\
 \tau^\mb{j} &= A^\mb{j}_{\mb{i} \ua} \om^{\mb{i}} \wedge \eta^\ua + B^\mb{j}_{\ua \ub} \eta^\ua \wedge \eta^\ub \\
  \tilde{\tau}^\ub &= \Cb{\mb{j} \mb{i}} \om^\mb{j} \wedge \om^\mb{i} + D^\ub_{\mb{j} \ua} \om^\mb{j} \wedge \eta^\ua + E^\ub_{\ua \uc} \eta^\ua \wedge \eta^\uc
\end{split}\]
with the understanding that each $\Cb{}, E^\ub$ and $B^\mb{j}$ are skew-symmetric.

Before diving into the main computation, we shall warm up by using our techniques to derive an integration by parts formula.  
\bgD{ehash}
For a differential operator $X$ on $\Sigma$ we define the horizontal adjoint $X^\#$ by
\[ \int f X h \; \Lambda_0  = \int h X^\# f  \; \Lambda_0 \]
for compactly supported functions $f,h$.
\enD

The key step to computing the horizontal adjoint of a vector field is the following computation on $\HF[,0]$:

 \bgE{divbit}
\begin{split}
d( \omj{j} \wedge \eta^* ) &= d \omj{j} \wedge \eta^* +(-1)^{k-1} \omj{j} \wedge d\eta^* \\
&= (-1)^{m-1} \om^{\mb{i}} \wedge \om^m_\mb{i} \wedge \omj{m,j} \wedge \eta^*\\ &
 \quad  + (-1)^{k+\ub} \omj{j} \wedge \tilde{\tau}^\ub \wedge \etj{\ub} \\
&= (-1)^{m-1} \om^0 \wedge \om^m_0 \wedge \omj{m,j} \wedge \eta^* + (-1)^{j+m} \om^m_j \wedge \omj{m} \wedge \eta^*\\
& \quad + (-1)^{k+\ub+j-1} 2\Cb{0j} \om^0 \wedge \om^* \wedge \etj{\ub}\\
&= (-1)^{m-1} \theta \wedge \om^m_0 \wedge \omj{m,j} \wedge \eta^* + (-1)^{j+m} \om^m_j \wedge \omj{m} \wedge \eta^*\\
& \quad + (-1)^{k+\ub+j-1}  2\Cb{0j}  \theta \wedge\om^* \wedge \etj{\ub}+(-1)^{j-1}2 a_\ub \Cb{0j} \om^* \wedge \eta^*.
\end{split}
\enE

\bgL{ehash}
For each $e_j$ with $j >0$ we have
\[ e_j^\# = -e_j -2a_\ub \Cb{0j} - \Gamma^m_{jm}\]
\enL

\pf
First note that \rfE{divbit} implies
\bgE{divbitS}
d(\omj{j} \wedge \eta^*)_{|\Sigma} = (-1)^{j-1} \left( 2a_\ub \Cb{0j} +\Gamma^m_{jm} \right) \Lambda_0. 
\enE
Thus
\begin{align*}
d( fh \omj{j} \wedge \eta^*)_{|\Sigma} &= (-1)^{j-1} ( f e_j h + h e_j f ) \Lambda_0 + (fh) d(\omj{j} \wedge \eta^*)_{|\Sigma}
\end{align*}
Thus 
\[\int  f e_j h \Lambda_0 = \int ( -he_j  f - 2fh a_\ub \Cb{0j} - fh \Gamma^m_{jm}  ) \Lambda_0 \]

\epf

Now we return to the derivation of a second variation formula. We begin by computing $\wh{F}^* d\theta = d \wh{F}^*\theta$ in two different ways and equating the results. Firstly 
\bgE{dt} \begin{split} 
d\theta &= d(\om^0 + a_\ub \eta^\ub) \\
&= \om^j \wedge \om_j^0 + \tau^0 + da_\ub \wedge \eta^\ub + a_\ub \tilde{\tau}^\ub\\
&= -\om^j \wedge \om^j_0  + A^0_{0 \ua} \om^0 \wedge \eta^\ua +A^0_{j \ua} \om^j \wedge \eta^\ua + da_\ub \wedge \eta^\ub\\
& \qquad + a_\ub 2\Cb{0j} \om^0 \wedge \om^j + a_\ub D^\ub_{0 \ua} \om^0 \wedge \eta^\ua +a_\ub D^{\ub}_{j \ua} \om^j \wedge \eta^\ua\\
& \qquad + \eta^\ua \wedge \eta^\ub \junk + \om^j \wedge \om^m \junk\\
&= -\om^j \wedge \om^j_0 +  \om^j \wedge \theta \big(-2a_\ub \Cb{0j} \big)  \\
& \quad + \om^j \wedge \eta^\ua \big(  A^0_{j \ua} +(e_j a_\ua) +a_\ub D^\ub_{j \ua} +2 a_\ua a_\ub \Cb{0j} \big)\\
& \qquad + \theta \wedge \eta^\ua \junk + \eta^\ua \wedge \eta^\ub \junk + \om^j \wedge \om^m \junk\\
\end{split}\enE
Thus 
\bgE{Fdt}
\begin{split}
\wh{F}^*d\theta &=  \om^j \wedge dt \big( -\uc^j_0 -2a_\ub \Cb{0j} \big)\\
& \quad + \om^j \wedge \eta^\ua \big(  -\Gm^j_{0 \ua} + A^0_{j \ua} +(e_j a_\ua) +a_\ub D^\ub_{j \ua}+2 a_\ua a_\ub \Cb{0j}  \big)\\
& \qquad + dt \wedge \eta^\ua \junk + \eta^\ua \wedge \eta^\ub \junk + \om^j \wedge \om^m \junk\\
\end{split}
\enE

But from the definitions we see that $\wh{F}^* \theta =\vfoc  dt$ so 
\[ d(\wh{F}^* \theta) = d\vfoc  \wedge dt = (e_j \vfoc ) \om^j \wedge dt + \eta^\ub \wedge dt \junk.\]

Comparing with \rfE{Fdt} thus yields
\bgE{gDer}
\begin{split}
e_j \vfoc  &= -\uc^j_0   -2\vfoc  a_\ub \Cb{0j} \\
0 &=-\Gm^j_{0 \ua} + A^0_{j \ua} +(e_j a_\ua) +a_\ub D^\ub_{j \ua} +2 a_\ua a_\ub \Cb{0j} 
\end{split}
\enE
Using metric compatibility of the connection thus yields the following useful identities
\bgE{guse}
\begin{split}
\uc^j_0 &=-e_j \vfoc   -2\vfoc  a_\ub \Cb{0j}\\
\Gamma^j_{0 \ua} &= A^0_{j \ua} + a_\ub D^\ub_{j \ua} +(e_j a_\ua) +2 a_\ua a_\ub \Cb{0j}  
\end{split}
\enE


Returning to the main computation. Recall that
\[ \Psi = (-1)^{j-1} \om^j_0 \wedge \omj{j} \wedge \eta^* \]
and so 
\bgE{tFsP}
\begin{split}
\partial_t \lrcorner  \wh{F}^* \Psi &= (-1)^{j-1} \uc^j_0 \wedge \omj{j} \wedge \eta^*.
 \end{split}
\enE
In particular, this implies that
\bgE{dtFsP}
\begin{split}
d( \partial_t \lrcorner \wh{F}^* \Psi) &= (-1)^{j-1} d\uc^j \wedge \omj{j} \wedge \eta^* + (-1)^{j-1} \uc^j_0 d (\omj{j} \wedge \eta^*) \\
d( \partial_t \lrcorner \wh{F}^* \Psi)_{|\Sigma} &=\big( e_j \uc^j_0  +\uc^j_0  \Gm^m_{jm}+  2 \uc^j_0 a_\ub \Cb{0j} \big) \Lambda_0\\
&= \big[ e_j (-e_j \vfoc   -2\vfoc a_\ub \Cb{0j} ) + (-e_j \vfoc   -2\vfoc a_\ub \Cb{0j} )( \Gm^m_{jm} +  2  a_\ub \Cb{0j} ) \big] \Lambda_0\\
&= \big[ -e_j^2 \vfoc   -4a_\ub \Cb{0j} e_j \vfoc  -\Gm^m_{jm} e_j \vfoc    -2\vfoc  e_j (a_\ub \Cb{0j})  \\
& \qquad -4\vfoc (a_\ub \Cb{0j})^2 -2\vfoc  \Gm^m_{jm}a_\ub \Cb{0j} \big] \Lambda_0.
\end{split}
\enE
Also if we define curvature $2$-forms by
\bgE{Curv}
d\om^\mb{i}_\mb{j} = \om^\mb{m}_\mb{j} \wedge \om_\mb{m}^\mb{i} + \Omega_\mb{j}^\mb{i}
\enE
then 
\bgE{dP}
\begin{split}
d\Psi &= (-1)^{j-1} d( \om^j_0 \wedge \omj{j} \wedge \eta^*)\\
&= (-1)^{j-1} \left( \om^m_0 \wedge \om_m^j + \Omega_0^j \right) \wedge \omj{j} \wedge \eta^* + (-1)^j \om_0^j \wedge d(\omj{j} \wedge \eta^*)\\
&=(-1)^{j-1} \left( \om^m_0 \wedge \om_m^j + \Omega_0^j \right) \wedge \omj{j} \wedge \eta^*\\
& \quad +  (-1)^{m+j} \theta \wedge \omega_0^j \wedge \om^m_0 \wedge  \omj{m,j} \wedge \eta^* + (-1)^{m}\omega^j_0 \wedge \om^m_j \wedge \omj{m} \wedge \eta^*\\
& \quad + (-1)^{k+\ub}  2\Cb{0j}  \theta \wedge \om_0^j \wedge \om^* \wedge \etj{\ub}- 2a_\ub \Cb{0j} \om_0^j \wedge \om^* \wedge \eta^*
\end{split}
\enE
Thus 
\bgE{tFdP}
\begin{split}
 \partial_t \lrcorner \wh{F}^*d\Psi  &= \Big[ \vfoc  \big( 2 \Omega^j_{00j} - \Gamma^j_{0m}\Gm^m_{0j}  +\Gamma^j_{0j}\Gamma^m_{0m}\\
 & \qquad -2\Cb{0j}\Gamma^j_{0\ub}\big)  -2a_\ub \Cb{0j}\uc^j_0   \Big]\om^* \wedge \eta^*
 \end{split}
\enE
When we restrict to $\Sigma$ we can use the fact that $H=-\Gamma^j_{0j} $ and $\vfo = \vfoc (\cdot,0)$ to see
\bgE{tFdPS}
\begin{split}
(\partial_t \lrcorner \wh{F}^*d\Psi )_{|\Sigma} &=  \vfo  \Big( \big(2\Omega^j_{00j} - \Gamma^j_{0m}\Gamma^m_{0j} +H^2 \big)\\
& \qquad  -2a_\ub \Cb{0j}(-e_j \vfo   -2\vfo a_\ub \Cb{0j})\\
& \qquad -2\vfo \Cb{0j}( A^0_{j \ub} + a_\uc  D^\uc_{j \ub} +(e_j a_\ub)  +2 a_\ub a_\ua \Cb[\ua]{0j}  ) \Big) \Lambda_0\\
&= e_j \vfo  \Big( 2a_\ub \Cb{0j}  \Big) \Lambda_0  \\
& \quad + \vfo  \Big(2\Omega^0_{00j} + \Gamma^0_{jm}\Gamma^0_{mj}  \\
& \qquad  -2\Cb{0j}A^0_{j \ub} - 2a_\uc D^{\uc}_{j \ub} \Cb{0j}- 2\Cb{0j} (e_j a_\ub) \Big) \Lambda_0
\end{split}
\enE

We now encode all this computation in the following lemma.

\bgL{Smooth}
Suppose $M$ is a vertically rigid subRiemannian manifold and $F$ is a noncharacteristic $C^{\infty;3}$ variation of $\Sigma\backslash C(\Sigma)$ with compactly supported horizontal variation function. Then
\bgE{Smooth}
\begin{split}
 \frac{d^2}{dt^2}\big |_{t=0} P_0( \Sigma_t ) 
 &= \int_{\Xi}   \left[  (\partial_t \vfoc ) H\Lambda \right]\big |_{t=0} + \int_{\Sigma} \left| \nabla^{0,\Sigma} \vfo  \right|^2 \Lambda  \\
& \quad + \int_{\Sigma} \vfo ^2 \big[ -\text{Ric}^{\nabla}(\nu,\nu) - \text{Tr}(\twoh^2) + H^2\\
 & \qquad   -\ap{\text{Tor}(\nu,e_j)}{T_\ub}{}\ap{\text{Tor}(e_j,T_\ub)}{\nu}{} \\
 & \qquad - \ap{\text{Tor}(e_j,T_\ub)}{N_v}{} \ap{\text{Tor}(\nu,e_j)}{T_\ub}{} \\
 & \qquad  - \Big( e_j \ap{\text{Tor}(\nu,e_j)}{N_v}{} -2\ap{\text{Tor}(\nu,e_j)}{N_v}{} \ap{ \nabla_{e_m} e_j}{e_m}{} \Big)\\
 & \qquad  - \ap{\text{Tor}(\nu,e_j)}{T_\ub}{} (e_j a_\ub) - \ap{\text{Tor}(\nu,e_j)}{N_v}{}^2\Big]\Lambda
 \end{split}\enE
\enL

\pf  As was shown \rfT[FV]{FirstVariation}
 \[ \frac{d}{dt} P_0( \Sigma_t ) =\int_{\Xi} \vfoc  F_t^* \Psi \]
 So
 \[\begin{split}
 \frac{d^2}{dt^2}_\big |{t=0} P_0( \Sigma_t ) &=\left[  \int_{\Xi} \Lie{\partial_t} (\vfoc  F_t^* \Psi)  \right]\big |_{t=0} \\
 &= \int_{\Xi}   \left[(\partial_t \vfoc ) H +  \vfo \Lie{\partial_t}(\wh{F}^* \Psi) \right]\big |_{t=0}.\\
 &= \int_{\Xi}   \left[  (\partial_t \vfoc ) H\Lambda +  \vfo \Lie{\partial_t}(\wh{F}^* \Psi) \right]\big |_{t=0}.\\
 \end{split}\]
Now as previously shown, locally 
 \bgE{Local}\begin{split} \left[ \Lie{\partial_t} \wh{F}^*\Psi \right]\bigg |_{t=0} &= \Big[ \sum\limits_{j=1}^k e_j^\# e_j \vfo     + \vfo  \big( 2\Omega^j_{00j} - \Gamma^0_{jm}\Gamma^0_{mj}+H^2  -2\Cb{0j}A^0_{j \ub} - 2a_\uc D^{\uc}_{j \ub} \Cb{0j}\\
 & \qquad  -2 e_j (a_\ub \Cb{0j}) - 2\Cb{0j} (e_j a_\ub) - 4 (a_\ub \Cb{0j} )^2\big) -2 \Gm^m_{jm}a_\ub \Cb{0j}\Big] \Lambda. \end{split}\enE
Converting to the invariant form given in the Lemma and integrating by parts once then completes the proof. 

\epf

\bgT{Main}
Suppose $M$ is a vertically rigid subRiemannian manifold, $\Sigma$ is a $C^2$ hypersurface and $F$ is a noncharacteristic $C^{1;3}$ variation of $\Sigma\backslash C(\Sigma)$. Then whenever either of the following holds
\begin{itemize}
\item $H=0$ on $\Sigma\backslash C(\Sigma)$
\item $H$ is constant on $C(\Sigma)$ and $F$ preserves $\int_\Xi F_t^* \mu$ for any smooth form with $\mu =dV$.
\end{itemize}
we have
\[\begin{split}
 \frac{d^2}{dt^2}\bigg |_{t=0} P_0( \Sigma_t ) 
 &= \int_{\Sigma} \left| \nabla^{0,\Sigma} \vfo  \right|^2 \Lambda  + \vfo ^2 \big[ -\text{Ric}^{\nabla}(\nu,\nu) - \text{Tr}(\twoh^2) \\
 & \qquad   -\ap{\text{Tor}(\nu,e_j)}{T_\ub}{}\ap{\text{Tor}(e_j,T_\ub)}{\nu}{} \\
 & \qquad - \ap{\text{Tor}(e_j,T_\ub)}{N_v}{} \ap{\text{Tor}(\nu,e_j)}{T_\ub}{} \\
 & \qquad  - \Big( e_j \ap{\text{Tor}(\nu,e_j)}{N_v}{} -2\ap{\text{Tor}(\nu,e_j)}{N_v}{} \ap{ \nabla_{e_m} e_j}{e_m}{} \Big)\\
 & \qquad  - \ap{\text{Tor}(\nu,e_j)}{T_\ub}{} (e_j a_\ub) - \ap{\text{Tor}(\nu,e_j)}{N_v}{}^2\Big]\Lambda
 \end{split}\]

 \enT
 
 \pf
 Since $F$ is supported away from the characteristic set, we have that 
 \[ \frac{d}{dt} \int_{\Xi} F_t^* \mu = \int_{\Xi} \Lie{\partial_t} F_t^* \mu = \int_{\Xi} \vfoc (\xi,t) \Lambda_t. \]
 Thus if $F$ preserves volume then $\int_{\Xi} \vfoc (\xi,t) \Lambda_t =0$. If $H$ is constant then differentiating yields
 \[ \int_{\Xi} (\partial_t \vfoc ) H \Lambda_0 + \int_{\Xi} \vfo ^2 H^2 \Lambda_0 =0. \]
 Therefore the effect of either condition is that the first term of \rfE{Smooth} cancels the $+\vfo ^2H^2$ term within the second integral. Thus the theorem is proved for $C^{\infty;3}$ variations.
 
 All that remains is to show that the result still holds with the restricted regularity conditions. The difficulty is that for the computations to hold, we must have $\vfo $ being $C^2$ on $\Xi$, whereas for a $C^{1;3}$ variation we can only guarantee that $\vfoc$ is continuous. However since $F_0$ itself is $C^2$ we see $\vfo$ is $C^1$. Fortunately, the right hand side of \rfE{Smooth}  requires only $C^1$ regularity in $\vfo $. All the other terms are in fact tensorial, so the restricted regularity will not cause problems.
 
 Now note that \[\wh{F}_t^* \Lambda = \lambda(\xi,t) d\xi^{1} \wedge \dots d\xi^{n}.\] Furthermore since $\Lambda$ is semibasic, we see by \rfC[BV]{rReg} that $\lambda$ is $C^{0;2}$. The second variation functional 
 \[ F \mapsto \frac{d^2}{dt^2}\bigg |_{t=0} P_0( F_t(\Xi) )\] is therefore continuous from $C^{1;3}$ variations to $\rn{}$.  By \rfL[BV]{Approx} we see that we can approximate $F$ by $C^{\infty;3}$ variations such that the restrictions to $t=0$ converge in $C^2$ to $F_0$. The second variation formula of \rfL{Smooth} holds for these approximations and the formula itself is continuous as a functional on $C^2$ embeddings.
  
\epf

\section{Examples}\setS{EX}

This second variation formula is hideously complicated in general so we shall attempt to illuminate it with some remarks and examples.

Firstly, recall that the horizontal second fundamental form is asymmetric but does have real valued entries. Thus its eigenvalues $\lambda_1,\dots \lambda_k$ are either real or come in conjugate pairs.  From elementary linear algebra we can then deduce
\bgE{EV}
\begin{split}
H &=  \text{Trace}(\twoh) =\sum_{j=1}^k \lambda_j\\
\text{Trace}(\twoh^2)&= \sum_{j=1}^k \lambda^2_j  = \sum_{j=1}^k \real{\lambda_j}^2 - \sum_{j=1}^k \imag{\lambda_j}^2
\end{split}
\enE
Now the imaginary parts of the eigenvalues reflect the asymmetry of \twoh which in turn reflects on the propensity of the  tangent horizontal vector fields to bracket generate the remaining directions. Reviewing the second variation formula of \rfT[SV]{Main} leads to the conclusion that a greater degree of bracket-generating causes greater stability in hypersurfaces.
  
\subsection{Strictly pseudoconvex pseudohermitian manifolds}

Recall a pseudohermitian manifold $(M,J,\eta)$ consists of:
\begin{itemize}
\item a $2n+1$-dimensional smooth manifold $M$
\item a non-vanishing $1$-form $\eta$ defining the horizontal distribution $V_0 = \text{ker} (\eta)$
\item a bundle map $J: V_0 \to V_0$ such that $J^2=-1$
\end{itemize}
with the integrability condition that the Niunhuis tensor (see \cite{Tanaka}) vanishes. The manifold is strictly pseudoconvex if the Levi metric
\[ g(X,Y) = d\eta(X,JY) + \eta(X) \eta(Y)\]
is positive definite. In this instance, a rigid vertical structure can be imposed by taking $T$ to be the Reeb vector field of $\eta$, i.e. $\eta(T)=1$ and $T \lrcorner d\eta=0$.

The second variation formula is then simplified by specializing the adapted connection to be the Tanaka-Webster connection (\cite{Tanaka} ,\cite{Webster}). Furthermore we require the horizontal frame to be $J$-graded, i.e. $Je_{2j}=e_{2j+1}$ and $J^* \om^{2j} =-\om^{2j+1}$, $j=0..n-1$. This implies that
\[ \tilde{\tau} = d\eta = \om^0 \wedge \om^1 + \om^2\wedge \om^3 + \dots \om^{n-2} \wedge \om^{n-1}.\]
Thus all the $D$ terms vanish, $\Cb[]{01}= \frac{1}{2}$  and $\Cb[]{0j}=0$ for $j>1$. 
Now
\[\begin{split} \Gm^m_{1m} &= \sum_{j=1}^{n-1} \left( \aip{ \nabla_{e_{2j}} J\nu }{e_{2j}}{} + \aip{ \nabla_{Je_{2j}} J\nu }{Je_{2j}}{} \right)\\
&= \sum_{j=1}^{n-1} \left( \aip{\nu }{ \nabla_{e_{2j}}  Je_{2j}}{} - \aip{  \nu }{\nabla_{Je_{2j}} e_{2j}}{} \right)\\
&= \sum_{j=1}^{n-1} \aip{[ e_{2j},Je_{2j}]}{\nu}{}\\
&= (n-1)a.
\end{split}\]
we the last line follows from the fact that each $[e_{2j},Je_{2j}]$ must be tangent to $\Sigma$, but $\eta([e_{2j},Je_{2j}])=-1$ together with the observation that $T-a\nu$ is tangent to $\Sigma$.

We can now use a combination of \rfT[SV]{Main} and its local expression \rfE[SV]{Local} to see
\[\begin{split}
 \frac{d^2}{dt^2}_{|t=0} P_0( \Sigma_t ) &=\int_{\Sigma} \Big[|\nabla^{0,\Sigma} \vfo |^2    + \vfo ^2 \big(-\text{Ric}^\nabla(\nu,\nu) -\text{Tr}(\twoh^2)  \\ & \qquad  -\ap{\text{Tor}(J\nu,T)}{\nu}{} 
  - 2 (J\nu)a  - na^2 \big) \Big] \Lambda 
\end{split}\]
If the pseudohermitian structure is normal (i.e. $\nabla_T X= [T,X]$, see \cite{Tanaka} for equivalent definitions and consequences) then the torsion term vanishes. For the case $n=1$ this example first appeared in \cite{CHMY}, although it should be noted that  their presentation of pseudohermitian manifolds causes $\Cb[]{01} =1$ instead.

A few particular examples are especially important in the literature:

\subsection{The Heisenberg Group} $\hn{n} = \rn[x,y]{2n} \times \rn[t]{}$ with the horizontal distribution spanned by
\[ X_j = \partial_{x^j} - \frac{1}{2} y^j \partial_t, \quad Y_j = \partial_{y^j} +\frac{1}{2} x^j \partial_t\]
is an example of both a (normal) strictly pseudoconvex pseudohermitian manifold (with $\eta = dt + \frac{1}{2} y^j dx^j - \frac{1}{2} x^j dy^j$, $JX_j =-Y_j$ and Reeb field $T=\partial_t$) and a Carnot Group. However, the curvature and the horizontal torsion both vanish identically so the second variation becomes
\[ \frac{d^2}{dt^2}_{|t=0} P_0( \Sigma_t ) =\int_{\Sigma} \Big[  |\nabla^{0,\Sigma} \vfo |^2    + \vfo ^2 \big( -\text{Tr}(\twoh^2)     - 2 (J\nu) a  -  na^2 \big) \Big] \Lambda 
\]
For $n=1$, this example was first shown by Danielli, Garofalo and Nhieu in \cite{DGN}. 

\subsection{The Rototranslation space}
Here $M = \rn{2} \times \sn{1}$ with horizontal distribution $V_0$ spanned by
\[ X_1 = \cos \theta \partial_x + \sin \theta \partial_y, \quad X_2=\partial_\theta.\]
A pseudohermitian structure can be imposed upon $M$, but we shall instead compute from first principles. Define \[T= [X_1,X_2] = \sin \theta \partial_x - \cos \theta \partial_y,\] so $[T,X_1]= 0$, $[T,X_2] =-X_1$. Thus $T$ represents a rigid vertical structure. A flat, adapted connection is created by defining \[\nabla_T X_1 =\nabla_T X_2 =0.\]  Direct computation shows that \[ \text{Tor}(X_1,X_2)= -T, \quad \text{Tor}(T,X_1)= 0, \quad \text{Tor}(T,X_2)=X_1\]
Now for a $C^2$ hypersurface $\Sigma$ we follow \cite{DGN:unstable} by writing the horizontal unit normal 
\[ \nu = \bt{p} X_1 + \bt{q} X_2, \quad e_1 = \bt{q}X_1 - \bt{p} X_2 \]
We can then compute that
\[ \ap{\text{Tor}(e_1,T)}{\nu}{}= -\bt{p}^2.\]
Applying  \rfT[SV]{Main}, we see that the second variation (of a minimal surface)  is given by \[\frac{d^2}{dt^2}_{|t=0} P_0( \Sigma_t ) =\int_{\Sigma} \Big[  (e_1\vfo )^2    + \vfo ^2 \big( -\bt{p}^2 - 2e_1 a  - a^2 \big) \Big] \Lambda 
\]
Now if we suppose $\Sigma$ is given as the level set of the defining function $\varphi$ and define $p = X_1 \varphi$, $q=X_2 \varphi$, $W = \nu \varphi$, $\omega= T\varphi$ and $\bt{\omega} = \omega/W$. Then since $a\nu -T$ is tangent to $\Sigma$ we have $a = \bt{\omega}$. Furthermore
\[\begin{split}
 e_1 a  &= \frac{1}{W} e_1 T \varphi - \frac{\bt{\omega}}{W} e_1 \nu\varphi\\
 &= \frac{1}{W} [e_1,T]\varphi - \frac{\bt{\omega}}{W} [e_1,\nu]\varphi\\
 &=  \frac{1}{W}  \left( -pT \bt{q} +q T \bt{p} - \bt{p}p \right) \\
 & \qquad -\frac{\bt{\omega}}{W} \left( q \nu \bt{p} -p \nu \bt{q} + \bt{q} [X_1 ,\nu] \varphi - \bt{p} [X_2, \nu] \varphi \right) \\
 &= \bt{q} T \bt{p} - \bt{p} T \bt{q}  -\bt{p}^2 - \bt{\omega} \left( \bt{q} \nu \bt{q} - \bt{p} \nu \bt{q} \right) \\
 & \qquad - \frac{\bt{\omega}}{W} \left(p \bt{q} X_1 \bt{p} + q \bt{q} X_1 \bt{q} + \bt{q}^2 \omega  - p \bt{p} X_2 \bt{p} -\bt{p} q X_2 \bt{q} + \bt{p}^2 \omega \right)\\
 &= \bt{q} T \bt{p} - \bt{p} T \bt{q}  -\bt{p}^2 + \bt{\omega} \left( \bt{p} \nu \bt{q}  -\bt{q} \nu \bt{q} \right) \\
& \qquad - \bt{\omega} \left(\frac{1}{2}\bt{q} X_1 (\bt{p}^2 +\bt{q}^2) -\frac{1}{2} \bt{p} X_2 (\bt{p}^2 +\bt{q}^2)  +\bt{\omega} \right)\\
&= \bt{q} T \bt{p} - \bt{p} T \bt{q}  -\bt{p}^2+ \bt{\omega} \left( \bt{p} \nu \bt{q}  -\bt{q} \nu \bt{q} \right) - \bt{\omega}^2
\end{split}\]
Thus, in this notation, the second variation formula becomes 
\[\frac{d^2}{dt^2}_{|t=0} P_0( \Sigma_t ) =\int_{\Sigma} \Big[  (e_1\vfo )^2    + \vfo ^2 \big( \bt{p}^2+ 2 (\bt{p} T \bt{q} - \bt{q} T \bt{p}) + 2\bt{\omega} ( \bt{q} \nu \bt{q} - \bt{p} \nu \bt{q} ) + \bt{\omega}^2 \big) \Big] \Lambda 
\]

\section{A compact, stable CMC surface in \hn{2}}\setS{HN}

In this section, we shall consider in detail the ``bubble sets'' in the five dimensional Heisenberg group \hn{2} and the connection to the isoperimetric problem. We define surface the $\Sigma \subset \hn{2}$ described as the completion of the double-graph 
 \[ t = \phi^{\pm}(r)=\pm \left( \frac{L^2 \pi}{8} -\frac{L^2}{4} \arctan \left( \frac{r}{\sqrt{L^2-r^2}} \right)+\frac{r}{4} \sqrt{L^2-r^2} \right), \quad 0 \leq r < L \]
 where $r^2 =\sum_{j=1}^{2} \left(x_j^2+y_j^2\right)$ and $L$ is a
 positive constant. This surface is easily seen to be $C^\infty$ away
 from its characteristic points at $(0,\pm L)$, and is $C^2$ but not
 $C^3$ over the characteristic locus \cite{dgn:iso}. It is also known
 to be the only rotationally invariant compact CMC surface in \hn{2},
 \cite{RR}. In \hn{1}, this bubble set was shown by Ritor\'e and
 Rosales to be the solution (in the category of $C^2$ surfaces) to the
 isoperimetric problem, \cite{RRas}.  In addition, Leonardi and Rigot \cite{LR}
 showed the bubble set in the first Heisenberg group to be the
 isoperimetric minimizer in the class of $C^2$ rotationally symmetric
 surfaces.  Very recently, Monti and Rickly \cite{MR} showed that the
 bubble sets in \hn{1} are isoperimetric minimizers in the class of
 convex surfaces.   In all Heisenberg groups,
 $\Sigma$ was shown by Danielli, Garofalo and Nhieu \cite{dgn:iso} to be the solution to the isoperimetric
 problem within the category of surfaces that can be described as
 double graphs over discs. 

We shall show that $\Sigma$ has constant mean curvature and is stable, in the sense that
 \[ \frac{d^2}{dt^2}\bigg |_{t=0} P_0(\Sigma_t) \geq 0 \]
under all volume preserving variations. Whilst not quite sufficient to establish that $\Sigma$ is a local minimum, this does at least provide further evidence that $\Sigma$ is a viable candidate for a solution to the isoperimetric problem in \hn{2}.

We start by defining the outward horizontal normal
\[ \nu = \Bt{p}^j X_j  +\Bt{q}^j Y_j  \] where \[W= \sqrt{ \sum \left( (X_j \phi)^2 +(Y_j \phi)^2 \right)} = \frac{Lr}{2\sqrt{L^2-r^2}}\] and 
\[ \Bt{p}^j= \frac{X_j \phi}{W} =  \frac{x^j}{L} - \frac{y^j\sqrt{L^2-r^2}}{Lr}, \quad \Bt{q}^j = \frac{Y_j \phi}{W} = \frac{y^j}{L} + \frac{x^j  \sqrt{L^2-r^2}}{Lr}.\]
Then $J\nu = \Bt{q}^j X_j - \Bt{p}^j Y_j$. We extend $\nu$, $J\nu$ to a horizontal frame by first defining
\begin{align*}
e &=  \frac{1}{r} \left( x^{2} X_1 -x^1 X_{1} - y^{2} Y_1 + y^1 Y_{2} \right) \\
Je &=\frac{1}{r} \left( y^1X_{2} - y^{2}X_1 + x^1 Y_{2}- x^{2}Y_1 \right).
\end{align*}
Now a straightforward, if brutal, computation shows that
 \[ \twoh = \begin{pmatrix} \frac{2}{L} & 0 &  0 \\ 0 & \frac{1}{L} & - \frac{\sqrt{L^2-r^2}}{Lr} \\   0 & \frac{\sqrt{L^2-r^2}}{Lr} & \frac{1}{L} \end{pmatrix} \]
  Thus the eigenvalues are
  \begin{align*}
  \kappa_0 &= \frac{2}{L} \\
  \kappa_1, \kappa_{2} &= \frac{1}{L} \pm i  \frac{\sqrt{L^2-r^2}}{Lr}\\
  \end{align*}
  So we can compute
  \begin{align*}
  H &= \frac{4}{L} \\
  \text{trace}\left(\twoh^2\right) &= \frac{6}{L^2} - \frac{2}{L^2r^2}(L^2-r^2)
  \end{align*}
  Now
  \[a= W^{-1}= \frac{2\sqrt{L^2-r^2}}{Lr}\]
  thus a simple computation shows
  \[ J\nu (a) = -\frac{2}{r^2}.\]
  Hence
  \begin{align*}
   - \text{trace}\left(\twoh^2\right) -2 J\nu(a) -2a^2 &= \frac{2}{L^2r^2} (L^2-r^2)  - \frac{6}{L^2} + \frac{4}{r^2} - 8 \frac{L^2-r^2}{L^2r^2}\\
   &= -\frac{2}{r^2}      \end{align*}
   So 
   \bgE{SV} \frac{d^2}{dt^2}_{|t=0} P_0(\Sigma_t) = \int_{\Sigma}\left( \left| \nabla^{0,\Sigma} \vfo \right|^2 - \frac{2}{r^2} \vfo^2 \right) \Lambda\enE
  at least for variations supported away from the characteristic locus. However, direct computation shows that
\[ |N_0| = \frac{Lr}{\sqrt{4L^2-4r^2+L^2r^2}} \]
and so $|N_0|$ is comparible to $r$ near $r=0$. For a variation
supported over the characteristic set, the horizontal variation
function can then be expressed as $\vfo =\frac{\vf}{|N_0|}$ and the
integrand of right hand side of \rfE{SV} is again integrable. Since
the characteristic locus is zero dimensional, the first variation formula holds without a boundary term even for variations supported over the characteristic locus. Reviewing the proof of the second variation formula reveals that the only place that required support away from the characteristic locus was the final integration by parts. If however, we cut out the characteristic locus by shrinking discs of radius $\e$ and apply \rfL[SV]{ehash}, we see that the boundary terms take the form $\vfo \nabla^{0,\Sigma} \vfo$ multiplied by a differential form bounded by the size of $\Lambda$. The boundary integrand blows up at a rate bounded by $\frac{1}{\e^2} = \frac{1}{\e^3} \cdot \e$. Due to the small size of the characteristic set, the boundaries of these discs are shrinking at rate $\e^3$ and so the boundary term is also negligible for characteristic variations. Hence \rfE{SV} holds for all variations.

To aid calculation, we now switch to polar coordinates $\rho_0 = \rho_0( r,\sigma)$ so  that \[ dV = \sqrt{1+ W^2 } r^{3} dr d\sigma\] and
  \[ \Lambda = \frac{ W }{\sqrt{1+W^2}} dV= \lambda(r) dr d\sigma\]
  where \[ \lambda(r) =  \frac{2r^{4}L}{\sqrt{L^2-r^2}}\]
  
 We note that away from the characteristic points, the level sets of the radial function $r$  foliate $\Sigma$ by $3$-dimensional spheres. Now $L^2(\sn{3},\cn{})$ can be orthogonally decomposed into (complex) homogeneous polynomials of bidegree $p,q$. Therefore on each foliating sphere we can decompose $\rho_0$ into spherical harmonics in $L^2(\sn{3})$,
  \[ \rho_0(r,\sigma) = \sum\limits_{p,q \geq 0} \rho^{p,q}_0 (r) \sigma_{p,q}. \]

\bgR{decomp}
 There is a natural pseudohermitian structure on the odd dimensional spheres and many natural subRiemannian operators have been studied using this orthogonal decomposition, see for example \cite{Folland:S}, \cite{Hladky2}. This technique can be generalized to study foliations by other compact pseudohermitian manifolds, \cite{Tanaka}, \cite{Hladky1}.
 \enR 
 
If we introduce the complex coordinates $z^j =x^j +iy^j$,   we can express $J\nu$ as,
  \[ J\nu = \sqrt{L^2-r^2} \partial_r + i \left( \bt{z}^j \partial_{\bt{z}^j} - z^j \partial_{z^j}\right)\] and so
  \[ J\nu (\rho_0) = \sum\limits_{p,q \geq 0} \left( \left[ \partial_r+i(q-p) \right]   \rho^{p,q}_0 \right) \sigma_{p,q}. \]
Now if $Z= \bt{z}^2 \partial_{z^1} -\bt{z}^1 \partial_{z^2}$ then
  \[ e =Z +\Bt{Z}, \quad Je =Z - \Bt{Z}.\]
 Now as a differential operator on $L^2(\sn{3})$, $Z^* =- \Bt{Z}$ so
  \begin{align*}
  \int_{\sn{3}} \left| e  \left( \rho_0(r,\sigma) \right)  \right|^2 &+ \left| Je \left( \rho_0(r,\sigma) \right) \right|^2 d\sigma \\
  &=\frac{1}{r^2}  \int_{\sn{3}}  -\rho_0 (Z +\Bt{Z})(Z+\Bt{Z})\rho_0 - \rho_0(\Bt{Z}-Z)(Z-\Bt{z}) \rho_0 d\sigma\\
  &=  \frac{1}{r^2} \int_{\sn{3}} -2 \rho_0 (Z\Bt{Z} + \Bt{Z}Z) \rho_0 d\sigma\\
  &= \frac{2}{r^2} \sum\limits_{p,q} \left[ p(q+1)+q(p+1)\right] \rho_0^{p,q}(r)
  \end{align*}
  Thus 
 \[ \int_{\Sigma} \sum_{m=1}^{n-1} \left(\left| e \rho_0 \right|^2 + \left| Je \rho_0 \right|^2 \right) \Lambda \geq \int_{\Sigma} \frac{2}{r^2} \left( \rho_0 - \rho_0^{0,0} \right)^2 \Lambda \]
 By spherical orthogonality, stability is proved if we can show
 \[ \int_{\Sigma} \left| J\nu \rho_0^{0,0} \right|^2 \Lambda \geq \int_{\Sigma} \frac{2}{r^2} \left(\rho_0^{0,0}\right) ^2 \Lambda \]
 Using the transformation $r=L \sin \theta$, $0 \leq \theta \leq \pi$, this reduces to showing
 \[ \int_0^\pi 2L^{2} \left( \partial_\theta \rho^{0,0}_0 \right)^2 \sin^{4} \theta d\theta \geq \int_0^\pi 4L^{2} \left( \rho^{0,0}_0 \right)^2 \sin^{2} \theta d\theta\]
 Now $\rho_0^{0,0}$ itself is generically singular at $\theta =0, \pi$ and indeed $\rho_0^{0,0}$ may not be in $L^2(0,\pi)$. But $h= \rho^{0,0}_0 \sin \theta$ is necessarily continuous up to $\theta=0,\pi$. It is also easy to see that $\partial_\theta h$ must be continuous and $L^2$ on $(0,\pi)$. 
 
 Recall that we are only interested in variations that preserve the volume, i.e. $\int_\Sigma \rho_0 \Lambda =0$. With these transformations this implies that $\int_0^\pi h \sin^{3} \theta \; d\theta =0$. 
 
 Thus stability for the geodesic ball in $\hn{n}$ is equivalent to the statement: for all $h \in C[0,\pi] \cap C^1(0,\pi)$ such that $h\upp \in L^2(0,\pi)$ and 
 \[ \int_0^\pi h \sin^{3} \theta \; d\theta =0\] the inequality
 \[ \int_0^\pi (h\upp \sin \theta - h \cos \theta)^2 \; d\theta \geq \int_0^\pi 2 h^2 \; d\theta\]
 holds.

The proof of this inequality is an elementary computation using a Fourier decomposition and is provided in the appendix. 


\bgR{higher}
It seems reasonable to conjecture that this result would also apply with $n>2$. However the required computations are more complicated due to the fact that for $n>2$, there is no useful global horizontal frame, even away from characteristic points. This is essentially the classical result that higher dimensional  spheres are not parallelizable. 
\enR

\appendix

\section{Size of the characteristic set}\setS{CS}

In this section we prove

\bgT{Main}
Suppose $X_1,X_2,\dots X_k$ are smooth globally defined vector fields on \rn{n} that bracket generate at every point. Then for any $C^2$ hypersurface $\Sigma$, the characteristic set
\[ C(\Sigma) = \{ p \in \Sigma \colon (X_j)\bigg|_{p} \in T_p \Sigma \text{ for all $1 \leq j \leq k$} \} \]
has Hausdorff dimension $\leq n-2$.
\enT 

The size and nature of the characteristic locus has been studied
widely \cite{Balogh,Magnani:cl,Derridj} in various contexts.  We include
a discussion here for completeness and because, to the best of our
knowledge, a complete argument for general sub-Riemannian spaces does
not appear in the current literature.

The proof is based on a series of technical lemmas.

\bgL{dimf=0}
Suppose $f$ is a $C^1$ function on \rn{m}. Then the set
\[ V =\{ p \colon f(p)=0, \; df\bigg |_{p} \ne 0 \} \] has Hausdorff dimension $\leq m-1$.
\enL

\pf
Fix $p \in V$. Then since the set $K= \{  p \colon f(p)=0, \; df\big |_{p} = 0\}$ is closed, we can find an open set $U$ containing $p$ such that $f_{|U}$ is a $C^1$ submersion from $U$ into \rn{}. The constant rank theorem implies that $U \cap V$ is a closed embedded submanifold of $U$. The set $U \cap V$ thus has Hausdorff dimension $m-1$ as a subset of $U$ (and hence as a subset of $\rn{m}$.) 

Therefore we can cover $V$ by open sets $U_\ua$ such that each $V \cap U_\ua$ has dimension $\leq m-1$. Since every subset of \rn{m} is second countable we can find a countable subcover by the Lindel\"of theorem. Thus we can express $V$ as a countable union of sets of dimension $\leq m-1$, which is sufficient to prove the result.
 
\epf

Our next lemma is a refinement of a result due to Derridj, Lemma 1 in \cite{Derridj}.

\bgL{dimBracket}
Suppose $\Sigma$ is a $C^2$ hypersurface in $\rn{n}$ and $X$ and $Y$ are smooth vector fields. Then the set
\[ V = \{ p \in \Sigma\colon X_p, Y_p \in T_p\Sigma, [X,Y]_p \notin T_p \Sigma\}\]
has Hausdorff dimension $\leq n-2$.
\enL

\pf
Locally we can introduce $C^2$ slice coordinates $(y,x^1,\dots ,x^{n-1})$ so that $\Sigma =\{ y=0\}$. Rewrite $X$ and $Y$ in these coordinates as
\begin{align*}
X &= a \partial_y  + a^i \partial_{x^i}\\
Y &= b \partial_y + b^i \partial_{x^i}.
\end{align*} 
Then 
\[ [X,Y] =  \left( a \pd{b}{y} - b\pd{a}{y} + a^i \pd{b}{x^i} - b^i \pd{a}{x^i} \right) \partial_y  \quad \text{mod } \partial_{x^1}, \dots \partial_{x^{n-1}}\]
The condition that $p=(0,p\upp) \in V$ is therefore equivalent to 
\[\begin{cases}
 & a(p)=0\\ & b(p)=0 \\
 & \left(a^i \pd{b}{x^i} - b^i \pd{a}{x^i}\right)(p)  \ne 0 
\end{cases}.\]
Set $a\upp = a_{| \Sigma}$, $b\upp=b_{|\Sigma}$. The portion of $V$ lying inside the slice coordinate chart must be contained in  $\{ a\upp =0, da\upp \ne 0\} \cup \{ b\upp =0, db\upp \ne 0\}$. Since $a\upp$ and $b\upp$ are (at least) $C^1$ functions on \rn{n-1}, the result now follows from \rfL{dimf=0}.
  
\epf

\noindent \textbf{Proof of \rfT{Main}: } Generate the countable collection of all vector fields that can be bracket generated by $X_1,\dots X_k$ and enumerate them as
\[ X_1,\dots X_k, \dots ,X_\ua , \dots \]
with the first $k$ matching the original vector fields. Define
\[ E_{\ua \ub} = \{ p \in \Sigma\colon (X_\ua)_{|p},  (X_\ub)_{|p} \in T_p\Sigma, [X_\ua,X_\ub]_p \notin T_p \Sigma\}.\]
Thus $\{ E_{\ua \ub} \}$ is a countable collection of sets of Hausdorff dimension $\leq n-2$. But since the original vector fields bracket generate at every point, for every $p \in C(\Sigma)$ we must be able to find $X_\ua$ and $X_\ub$ such that $p \in E_{\ua \ub}$. Therefore $C(\Sigma)$ is contained in the countable union of sets of Hausdorff dimension $\leq n-2$ and so must have dimension $\leq n-2$ also.

\epf
 
\bgR{Sharp}
Without further restrictions on the vector field $X_j$ this result is sharp for all $n>k \geq 2$. To see this set
\begin{align*}
X_j & = \partial_{x^j} \quad \text{for $1 \leq j \leq k-1$}\\
X_k & = \partial_{x^k} + (x^1) \partial_{x^{k+1}} + (x^1)^2 \partial_{x^{k+2}} + \dots + (x^1)^{n-k} \partial_{x^n}.
\end{align*}
These vector fields bracket generate at step $n-k+1$ at all points of \rn{n}. The smooth surface $\Sigma = \{ x^n = (x^1)^2\}$ then has the property \[\{ x^n =x^1 =0 \} \subset C(\Sigma) .\]
Thus the Hausdorff dimension of $C(\Sigma)$ must be $\geq n-2$.
\enR

In the special case of the higher Heisenberg groups this theorem is decidedly non-sharp. It was shown by Balogh, \cite{Balogh}, that for the Heisenberg group of dimension $2n+1$ the characteristic set dimension is bounded by $n$ rather than $2n-1$. Balogh also showed that if the condition $C^2$ is relaxed to $C^{1,1}$ then the bound $<2n$ is actually sharp. 

The improved bounds for the higher Heisenberg groups was independently shown by Cheng-Hwang \cite{CHY} for graphs over the horizontal variables. Their technique had the advantage that it used only elementary linear algebra and generalized to graphs in pseudohermitian manifolds in natural coordinates. Here we present a new coordinate free version of this approach which can be used as a tool to study characteristic dimension in general equiregular subRiemannian structures.

\bgD{Hess} Given a collection of vector fields $\mathcal{X}=\{X_1, \dots X_n\}$ and a $C^2$ function $\phi$, we define the Hessian of $\phi$ at $p$ with respect to $\mathcal{X}$ by
\[ \hess{2}(\phi,p) = \left( X_j X_k \phi _{|p} \right)\] Additionally we define the  symmetric Hessian and skew-symmetric Hessian by
\[  \hess[+]{2} = \hess{2} + \left( \hess{2} \right)^{\top}, \quad \hess[-]{2} = \hess{2}- \left( \hess{2} \right)^{\top}. \]
Thus we note
\bgE{Hess} \hess[-]{2}(\phi,p) = \left( [X_j,X_k]_{|p} \phi \right).\enE
\enD

Now at any point $p$, we can find a non-degenerate constant matrix $P$ such that the skew-symmetric Hessian can be written
\[ \hess[-]{2} = P^{-1} \begin{pmatrix} J_\lambda & 0 \\ 0 & 0 \end{pmatrix} P\]
where $J_\lambda$ is the $2\lambda \times 2 \lambda$ matrix with $\lambda$ copies of $\begin{pmatrix} 0 & 1\\ -1 & 0\end{pmatrix}$ along the leading diagonal and zeros everywhere else. Thus we immediately obtain that 
\[ \text{rank} \left(\hess[-]{2}\right) = 2\lambda, \quad \text{rank}\left( \hess{2}  \right)\geq \lambda.\]

As a basic illustration of the use of these Hessian we present the following lemma, which essentially first appeared in \cite{CHY}.

\bgL{pseudo}
Suppose $(M,\eta,J)$ is a $2m+1$ dimensional pseudohermitian manifold such that the Levi form 
\[ (X,Y) \mapsto d\eta (X,JY)=- \eta[X,JY], \qquad \text{$X$, $Y$ horizontal} \]
has signature $(p,n)$ with $ p+n \geq 2k$ everywhere. Then is $\Sigma$ is any $C^2$ hypersurface, the characteristic set $C(\Sigma)$ has Hausdorff dimension $\leq 2m-k$.
\enL

\pf
Choose $p \in C(\Sigma)$ and let $\phi$ be a $C^2$ defining function for $\Sigma$ in a neighbourhood of $p$. Next choose $\mathcal{X}= \{ X_1, \dots X_{2m} \}$ with $X_{m+j}= J X_m$ a frame for the horizontal distribution near $p$. Now $T \phi$ cannot vanish at $p$ as otherwise $d\phi_{|p}=0$. Since the Levi form has total signature bounded below by $2k$, from \rfE{Hess} have that $\text{rank} \left(\hess[-]{2}(\phi,p) \right) \geq 2k$. Thus $\text{rank}\left(\hess{2}(\phi,p) \right)  \geq k$.

Define a function $F\colon M \to \rn{2m+1}$ by 
\[ F = \begin{pmatrix} \phi \\ X_j \phi \end{pmatrix}. \] If we extend $\mathcal{X}$ to by a vector field $T$ to a frame for $M$ near $p$ we see
\[ DF_p = \begin{pmatrix}  0 & T\phi \\ \hess{2}(\phi,p) & * \end{pmatrix}\] thus $DF_p$ has rank $\geq k+1$. But near $p$, $C(\Sigma) = F^{-1}(0)$ so the intersection of $C(\Sigma)$ with a neighbourhood of $p$  is contained in an embedded submanifold of dimension $\leq 2m-k$.
 
\epf

This technique can be extended to equiregular subRiemannian manifolds of higher step or otherwise more complicated vertical structures, but the generically the derived bounds on characteristic dimension are no better than the general result of \rfT{Main}.

Suppose $M$ is an $n$-dimensional equiregular subRiemannian manifold $M$ 
with $\mathcal{X}$ a smooth local frame for the horizontal distribution. Then we can produce smooth frames $\mathcal{X}_{(1)}, \mathcal{X}_{(2)}, \dots$ consisting of vector fields produced from $\mathcal{X}$ by $1,2, \dots$ or less commutations respectibely. Then for any $C^2$ surface $\Sigma$ and any point $p \in C(\Sigma)$ we again study the skew-symmetric Hessian $\hess[-]{2}(\phi,p)$. If this Hessian does not vanish at $p$ we can deduce that $C(\Sigma)$ is locally contained in an embedded submanifold of dimension $\leq n-2$. If the $\hess[-]{2}$ does vanish we can immediately deduce that $p$ is actually a characteristic point for the distribution $\mathcal{X}_{(1)}$. We then iterate this argument. If $\mathcal{X}$ is bracket-generating, this must terminate and we have rederived the result of \rfT{Main}.

If the step size of the subRiemannian structure is greater than $1$, then in particular this argument suggests that generically we cannot expect any improvement on the bound $n-2$. This would not be surprising as the condition that the dimension of the the hypersurface equaling the dimension of the horizontal distribution might be expected to yield a richer theory than the general case. That said, there are examples where this technique can produce improved bounds.

\bgX{HH2}
Consider $M = \hn{\hn{2}}$ the $21$-dimensional manifold constructed as follows: 
\[ X_{j,k}= \partial_{x_{j,k}} - \frac{1}{2} y_{j,k} U_j, \quad Y_{j,k}= \partial_{y_{j,k}} + \frac{1}{2} x_{j,k} U_j , \quad j,k =1..4\] 
with the $U_j$ the horizontal generators of an independent copy of $\hn{2}$. Thus $M$ consists of $4$ copies of $\hn{2}$ each yielding an element of another $\hn{2}$ as its characteristic field. $M$ is then a step $2$ Carnot group with codimension $5$ horizontal distribution. 

If $\Sigma$ is a $C^2$ hypersurface with defining function $\phi$ then for any $p \in C(\Sigma)$, either $\hess[-]{2}(\phi,p)$ vanishes identically or has rank $\geq 4$. However if $\hess[-]{2}$ vanishes identically then $p$ is a characteristic point for $\mathcal{X}_{(1)}$. But this higher level characteristic set is contained in a submanifold of dimension $\leq 21- 3$ by an identical argument.
\enX

\section{Technical Lemma}\setS{TN}

This section is devoted to proving the following technical lemma.

\bgL{Tech}
Suppose that $\rho \in C[0,\pi] \cap C^1(0,\pi)$ and $\pr{\rho} \in L^2(0,\pi)$. If 
\[ \int_0^\pi \rho \sin^3 \theta d\theta =0\]
then \[ \int_0^\pi \left(\pr{\rho}  \sin \theta - \rho \cos \theta \right)^2 d\theta \geq 2\int_0^\pi \rho^2 d\theta.\]
\enL

\pf 
We first note that
\[ \left(\pr{\rho}  \sin \theta - \rho \cos \theta \right)^2 = \left( [\rho \sin\theta]\upp \right)^2 - 4\rho \rho\upp \sin \theta \cos \theta\]
Now by an integration by parts
\begin{align*}
 \int_0^\pi 4\rho \rho\upp \sin \theta \cos \theta d\theta &= - \int_0^\pi 2\rho^2 \left( \cos^2 \theta - \sin^2 \theta \right) d\theta\\
 &= -\int_0^\pi 2 \rho^2 d\theta +4 \int_0^\pi \left(\rho \sin \theta\right)^2 d\theta
\end{align*}
Thus
\bgE{trans}\int_0^\pi \left(\pr{\rho}  \sin \theta - \rho \cos \theta \right)^2 -2 \rho^2 d\theta = \int_0^{\pi} \left([\rho \sin \theta]\upp\right)^2 - 4 \left(\rho \sin \theta\right)^2 d\theta.\enE
The lemma is proved if we can prove that this integral is positive. We proceed by splitting $\rho \sin \theta$ into its Fourier decomposition on $(0,\pi)$, i.e. as an element of $W^1(0,\pi)$
\[ \rho \sin \theta = \frac{a_0}{2} + \sum\limits_{k \geq 1} \left( a_k \cos(2k \theta) + b_k \sin(2k\theta) \right).\]
Therefore the integrals in \rfE{trans} equal
\bgE{sum}
 \frac{\pi}{2} \left(- \frac{1}{2} a_0^2+ \sum\limits_{k \geq 1}  (4k^2-4) (a_k^2 +b_k)^2  \right).
 \enE
Now the problem is the constant term $a_0$, however we have yet to use the fact that
\[ \int_0^\pi (\rho \sin \theta) \sin^2\theta d\theta = 0.\]
Writing $\sin^2 \theta = \frac{1}{2} (1-\sin(2\theta))$ we see that $a_0 =a_2$. Thus using the term $k=2$ we see
\[ -\frac{a_0^2}{2} + (16 -4)(a_2^2+b_2)^2 \geq 0. \]
The integrals of \rfE{trans} are therefore positive and the lemma is proved.

\epf

\bgC{tech}
The only function $\rho$ satisfying the conditions of \rfL{Tech} such that the inequality is an equality is $\rho= a_1 \cos \theta$.
\enC

\pf
For $k \ne 1$, the inequality at each level of the sum in \rfE{sum} is strict. Thus $\rho$ must have $a_k = b_k =0$ for $k \ne 0$. But if $b_1 \ne 0$ then $\rho \sin \theta$ is not orthogonal to $\sin^2 \theta$.

\epf
\bibliographystyle{plain}
\bibliography{References}

\end{document}